\documentclass[amsthm]{elsart}

\usepackage{yjsco}
\usepackage{natbib}

\usepackage[table]{xcolor}
\usepackage{amsmath}
\usepackage{amssymb}
\usepackage{amsthm}
\usepackage{dsfont}
\usepackage{mathtools}
\usepackage{threeparttable}
\usepackage{hyperref}
\usepackage{placeins}
\usepackage{pgfplots}
\usepackage{pgfplotstable}
\usepackage{subcaption}
\usepackage{wrapfig}
\usepackage[title]{appendix}
\usepackage{algorithm}
\usepackage{algpseudocode}
\usepackage{graphicx}

\usepackage{amssymb,amsmath}
\usepackage{pgfplots}
\usepackage{pgfplotstable}

\newcommand{\term}[1]{\textbf{#1}}
\newcommand{\name}[1]{#1}
\newcommand{\techy}[1]{\textsf{#1}}
\newcommand{\op}[1]{\ensuremath{\operatorname{#1}}}
\newcommand{\zmat}[2]{\ensuremath{M_{#1 \times #2}(\Z)}}

\newcommand{\boldx}{\boldsymbol{x}}
\newcommand{\boldy}{\boldsymbol{y}}
\newcommand{\boldz}{\boldsymbol{z}}
\newcommand{\bolda}{\boldsymbol{a}}
\newcommand{\boldb}{\boldsymbol{b}}

\newcommand{\boldp}{\boldsymbol{p}}
\newcommand{\boldt}{\boldsymbol{t}}
\newcommand{\boldalpha}{\boldsymbol{\alpha}}
\newcommand{\boldbeta}{\boldsymbol{\beta}}

\newcommand{\boldone}{\boldsymbol{1}}
\newcommand{\boldzero}{\boldsymbol{0}}
\newcommand{\N}{\mathbb{N}}
\newcommand{\Z}{\mathbb{Z}}
\newcommand{\R}{\mathbb{R}}
\newcommand{\C}{\mathbb{C}}
\newcommand{\Cstar}{\mathbb{C}^\ast}
\newcommand{\V}{\mathcal{V}^\ast}

\newcommand{\calD}{\ensuremath{ \mathcal{D} }}

\newcommand{\hatS}{\ensuremath{ \hat{\mathcal{S}} }}

\newcommand{\nlist}[1]{\ensuremath{#1_{1},\dots,#1_{n}}}

\newcommand{\Vol}{\op{Vol}}
\newcommand{\nvol}{\op{NVol}}
\newcommand{\conv}{\op{conv}}
\newcommand{\known}{\cellcolor{lightgray}}

\newtheorem{theorem}{Theorem}
\newtheorem{proposition}{Proposition}
\newtheorem{conjecture}{Conjecture}

\theoremstyle{definition}
\newtheorem{definition}{Definition}
\newtheorem{remark}{Remark}
\newtheorem{example}{Example}

\begin{document}
\begin{flushright}
ADP-15-2/T904
\end{flushright}

\begin{frontmatter}

\title{
    Parallel degree computation for solution space of binomial systems
    with an application to the master space of $\mathcal{N}=1$ gauge theories
}

\author{Tianran Chen}
\address{Department of Mathematics, Michigan State University, East Lansing, MI, USA}
\ead{chentia1@msu.edu}

\author{Dhagash Mehta}
\address{Department of Applied and Computational Mathematics and Statistics, University of Notre Dame,
Notre Dame, IN 46545, USA}
\address{Centre for the Subatomic Structure of Matter, 
Department of Physics, School of Physical Sciences, University of Adelaide, Adelaide, South Australia 5005, Australia}
\ead{dmehta@nd.edu}



\begin{abstract}
    The problem of solving a system of polynomial equations is one of the 
    most fundamental problems in applied mathematics. 
    Among them, the problem of solving a system of binomial equations form
    a important subclass for which specialized techniques exist.
    For both theoretic and applied purposes, the degree of the solution set
    of a system of binomial equations often plays an important role in understanding
    the geometric structure of the solution set. Its computation, however, 
    is computationally intensive.
    This paper proposes a specialized parallel algorithm for computing the degree
    on GPUs that takes advantage of the massively parallel nature of GPU devices.
    The preliminary implementation shows remarkable efficiency and scalability
    when compared to the closest CPU-based counterpart.
    Applied to the ``master space problem of $\mathcal{N}=1$ gauge theories''
    the GPU-based implementation achieves nearly 30 fold speedup over its
    CPU-only counterpart enabling the discovery of previously unknown results.
    Equally important to note is the far superior scalability:
    with merely 3 GPU devices on a single workstation, the GPU-based implementation
    shows better performance, on certain problems, than a small cluster totaling
    100 CPU cores.
\end{abstract}


\end{frontmatter}

\maketitle

\section{Introduction}

The problem of solving a system of polynomial equations
is one of the most fundamental problems in applied mathematics and science.
Among them, the problem of solving a system of binomial equations
is of special interest for they appear naturally in many applications, 
and specialized and much more efficient algorithm exists (e.g. \cite{kahle_decompositions_2010}).
In many applications, only the solutions of a system of binomial equations
for which no variable is zero are needed.
Such solutions are known as the $\Cstar$-solution set
and will be the focus of this article.

In this article, we propose a parallel algorithm for computing the 
degree of a $\Cstar$-solution set of a system of binomial equations.
This algorithm is specially designed for GPU (graphics processing unit) devices
by taking advantage of the massively parallel nature of GPUs.
When applied to a binomial system coming from particle physics, 
called the master space of $\mathcal{N}=1$ gauge theories,
this algorithm is able to produce previously unknown results.
Furthermore, the experimental implementation for GPU built on top of the 
\techy{CUDA} framework has already shown promising results.
Remarkably, with multiple GPU devices (on the same computer),
the GPU based implementation exhibits much better performance,
in many cases, than small to medium sized computer clusters.

This article is structured as follows:
First, necessary notations and concepts are introduced.
In particular, we shall review basic geometric properties of 
the $\Cstar$-solution set defined by a binomial system.
Then the algorithm for computing the Smith Normal Form of an integer matrix
is reviewed in \S \ref{sec:computation}, as it is an important tool necessary
in understanding the structure of the $\Cstar$-solution set of a binomial system.
The core of this article is \S \ref{sec:degree} where a highly scalable parallel 
algorithm for computing the degree of the $\Cstar$-solution set of 
a system of binomial equations is presented.
A natural by-product of the degree computation is a series of homotopy constructions
that can be used to compute the ``witness sets'' of components of the 
$\Cstar$-solution set of a binomial system, which is an important and ubiquitous 
construction in numerical algebraic geometry.
This process is explained in \S \ref{sec:witness}.
The problem of studying the master space of $\mathcal{N}=1$ gauge theories, arising from string theory, 
is briefly reviewed in \S \ref{sec:mspace},
and we show the interesting and previously unknown results obtained by applying 
the parallel algorithm for solving systems of binomial equations
and computing the degree of the solution set to the master space problem.

\section{Laurent binomial systems and its solution set} 
\label{sec:binomial}

First, we shall introduce necessary concepts and notations.
For positive integers $m$ and $n$,
let $\zmat{n}{m}$ denote the set of all $n \times m$ integer matrices.
A square integer matrix is said to be \term{unimodular} if its determinant is $\pm 1$.
Note that such a matrix $A \in \zmat{n}{n}$ has a unique inverse 
$A^{-1} = \frac{1}{\det A} \op{adj} A$ that is also in $\zmat{n}{n}$,
where $\op{adj} A$ is the adjoint matrix of $A$.
The $n \times n$ identity matrix in $\zmat{n}{n}$ is denoted by $I_n$.

Even though the main application considered in this article are binomial systems,
its theory is more naturally developed in the context of more general 
``Laurent binomial systems'' where negative exponents are allowed.
For variables $\boldx = (x_1, \dots, x_n)$, a \term{Laurent monomial} in $\boldx$ 
is an expression of the form $x_1^{\alpha_1} \cdots x_n^{\alpha_n}$
where $\alpha_1, \dots, \alpha_n$ are integers (which may be zero or negative). 
For convenience, we shall write $\boldalpha = (\alpha_1,\dots,\alpha_n)^\top \in \Z^n$
and use the ``vector exponent'' notation
\[ 
    \boldx^{\boldalpha} = 
    (x_1, \dots, x_n)^{ 
        \left(
            \begin{smallmatrix} \alpha_1 \\ \vdots \\ \alpha_n \end{smallmatrix} 
        \right)
    } = 
    x_1^{\alpha_1} \, \cdots \, x_n^{\alpha_n}
\]
to denote a Laurent monomial.
Similarly, for an integer matrix $A \in M_{n \times m}(\Z)$
with columns $\boldalpha^{(1)},\dots,\boldalpha^{(m)} \in \Z^n$,
the ``matrix exponent'' notation will be used for an $m$-tuple of Laurent monomials:
\begin{equation}
    \label{equ:matrix-power}
    \boldx^A \; = \;
    \boldx^{ \begin{pmatrix} \boldalpha^{(1)} & \cdots & \boldalpha^{(m)} \end{pmatrix} } \; := \;
    ( \boldx^{\boldalpha^{(1)}}, \dots, \boldx^{\boldalpha^{(m)}} ) .
\end{equation}
This notation is particularly convenient since the familiar identities 
$\boldx^{I_n} = \boldx$ and $(\boldx^A)^B = \boldx^{AB}$ still hold.
Since the exponents here may be negative,
it is only meaningful to consider the function $\boldx \mapsto \boldx^A$
when we restrict each $x_i$ to be nonzero. In particular, throughout this article, 
we shall let $x_i \in \C^\ast = \C \setminus \{0\}$ for each $i=1,\dots,n$.
In this case, each matrix $A \in \zmat{n}{m}$ induces a function
from $(\C^\ast)^n$ to $(\C^\ast)^m$ given by $\boldx \mapsto \boldx^A$.
Of particular importance is the function induced by a unimodular matrix $A \in \zmat{n}{n}$ 
since $A^{-1}$ is also in $\zmat{n}{n}$, and hence functions 
$\boldx \mapsto \boldx^A$ and $\boldx \mapsto \boldx^{A^{-1}}$ are the inverses of each other
($(\boldx^A)^{A^{-1}} = \boldx^{A \, A^{-1}} = \boldx^{I_n} = \boldx$).

A \term{Laurent binomial} is 
an expression of the form
$c_1 \boldx^{\boldalpha} + c_2 \boldx^{\boldbeta}$
for some $c_1, c_2 \in \C^\ast$ and $\boldalpha, \boldbeta \in \Z^n$.
This article focuses on the properties of the solution set of systems of Laurent binomials equations,
or simply \term{Laurent binomial systems}, over $(\C^\ast)^n$.
Stated formally, given exponent vectors 
$\boldalpha^{(1)}, \dots, \boldalpha^{(m)}, \allowbreak \boldbeta^{(1)}, \dots, \boldbeta^{(m)} \in \Z^n$ 
and the coefficients $c_{i,j} \in \C^\ast$,
the goal is to describe the set of all $\boldx \in (\C^\ast)^n$ 
that satisfies the system of equations
\[ 
    \begin{cases}
    	c_{1,1} \boldx^{\boldalpha^{(1)}} + c_{1,2} \boldx^{\boldbeta^{(1)}} &= 0 \\
    	&\vdots \\
    	c_{m,1} \boldx^{\boldalpha^{(m)}} + c_{m,2} \boldx^{\boldbeta^{(m)}} &= 0
    \end{cases}.
\]
Since only the solutions in $(\C^\ast)^n$ are concerned, 
this system is clearly equivalent to
\[ 
    (
	\boldx^{\boldalpha^{(1)} - \boldbeta^{(1)}},
	\dots,
	\boldx^{\boldalpha^{(m)} - \boldbeta^{(m)}} 
    ) \; = \;
    (-c_{1,2}/c_{1,1}, \dots, -c_{m,2}/c_{m,1}).
\]
With the more compact ``matrix exponent'' notation in \eqref{equ:matrix-power},
this system can simply be written as
\begin{equation} \label{equ:standard-form}
    \boldx^A = \boldb \qquad \text{ or equivalently } \qquad \boldx^A - \boldb = \boldzero
\end{equation}
where the integer matrix $A \in M_{n \times m}(\Z)$,
having columns $\boldalpha^{(1)}-\boldbeta^{(1)},\dots,\boldalpha^{(m)}-\boldbeta^{(m)}$,
represents the exponents appeared in the Laurent monomials and the vector 
$\boldb = (-c_{1,2}/c_{1,1}, \linebreak[1] \dots, \linebreak[1] -c_{m,2}/c_{m,1})^\top \in (\C^\ast)^m$
collects all the coefficients.
The solution set of \eqref{equ:standard-form} over $(\C^\ast)^n$ shall be denoted by
\begin{equation} \label{equ:soln-set}
    \V(\boldx^A-\boldb) = \{ \boldx \in (\C^\ast)^n \mid \boldx^A - \boldb = \boldzero \}.
\end{equation}
The goal of this article is to present efficient parallel algorithms for computing 
the structural properties of the set $\V(\boldx^A-\boldb)$:
its dimension, number of components, global parametrizations, and, most importantly, degree.
We shall first briefly review some basic facts about the $\C^\ast$-solution set 
of a Laurent binomial system.
A more detailed summary can be found in the article \cite{chen_solutions_2014} 
by the first author and \name{Tien-Yien Li}.
In depth theoretical discussions can be found in standard references such as
\cite{cox_toric_2011,eisenbud_binomial_1996,fulton_introduction_1993,miller_combinatorial_2005,sturmfels_equations_1997}.
Certain computational aspects have been studied in \cite{kahle_decompositions_2010,kahle2014decompositions}.

An important tool in understanding the structure of $\V(\boldx^A-\boldb)$
is the \emph{Smith Normal Form} of the exponent matrix $A \in \zmat{n}{m}$:
there are unimodular square matrices 
$P \in \zmat{n}{n}$ and $Q \in \zmat{m}{m}$ such that
\begin{equation}
    \label{equ:smith}
    P \, A \,\, Q \; = \; 
    \begin{pmatrix}
    	d_1              &        & \phantom{\ddots} &   &        &  \\
    	\phantom{\ddots} & \ddots & \phantom{\ddots} &   &        &  \\
    	                 &        & d_r              & \phantom{\ddots}  &        &  \\
    	                 &        &                  & 0 &        &  \\
    	                 &        &                  &   & \ddots & \phantom{\ddots} \\
    	                 &        &                  &   &        & 0
    \end{pmatrix}
    \in \zmat{n}{m}
\end{equation}
with nonzero integers $d_1 \, | \, d_2 \, | \, \cdots \, | \, d_r$ for $r = \op{rank} A$,
unique up to the signs. 
Here, $a \,|\,b$ means $a$ divides $b$ as usual.
This decomposition of the matrix $A$ provides important topological information 
about $\V(\boldx^A - \boldb) \subset (\Cstar)^n$ 
summarized in the following proposition:

\begin{proposition}[Topological description \cite{eisenbud_binomial_1996}] \label{pro:topological}
    If $\V(\boldx^A - \boldb)$ in $(\C^\ast)^n$ is not empty,
    then it consists of a finite number of connected components. Furthermore,
    \begin{enumerate}
    	\item the number of components is exactly $\left| \prod_{j = 1}^r d_j \right|$.
    	\item each solution component has codimension equal to $\op{rank} A = r$.
    \end{enumerate}
\end{proposition}

This description can be strengthened significantly.
Here we shall briefly outline the derivation of the stronger description 
of $\V(\boldx^A - \boldb)$ as it provides the important data that form the starting point
of the degree computation to be discussed in \S \ref{sec:degree}.
For $P$ and $Q$ in the Smith Normal Form of $A$ in (\ref{equ:smith}), 
let $P_r \in M_{r \times n}(\Z)$ and $P_0 \in M_{(n-r) \times n}(\Z)$ 
be the top $r$ rows and the remaining $n - r$ rows of $P$ respectively.
Similarly, let $Q_r \in M_{m \times r}(\Z)$ and $Q_0 \in M_{m \times (m-r)}(\Z)$
be the left $r$ columns and the remaining $m-r$ columns of $Q$ respectively.
With these notations, the Smith Normal Form of \eqref{equ:smith} of $A$ can be written as
\begin{equation}
    \label{equ:rank-decomp}
    \begin{pmatrix}
        P_r \\ P_0
    \end{pmatrix}
    A
    \begin{pmatrix}
        Q_r & Q_0
    \end{pmatrix}
    =
    \begin{pmatrix}
        D         & \boldzero \\
        \boldzero & \boldzero
    \end{pmatrix}
\end{equation}
with $D = \op{diag}(d_1,\dots,d_r) \in M_{r \times r}(\Z)$
and $\boldzero$'s representing zero block matrices of appropriate sizes.
With this we can transform the binomial system $\boldx^A = \boldb$
into a form from which more detailed information can be easily extracted.

Since $P$ and $Q$ are both unimodular 
the maps $\boldz \mapsto \boldz^P$ and $\boldy \mapsto \boldy^Q$ are 
both bijections on $(\C^\ast)^n$ and $(\C^\ast)^m$ respectively.
Therefore, considering the solution set in $(\C^\ast)^n$, the original system $\boldx^A = \boldb$ is equivalent to
$(\boldx^A)^Q \; = \; \boldx^{A\,Q} \; = \; \boldb^Q$.
Similarly, the solution sets remain equivalent after the change of variables 
$\boldx = \boldz^P$, which produces
\[
    (\boldz^P)^{A\,Q} \; = \; 
    \boldz^{PA\,Q} \; = 
    \boldz^{
	\begin{psmallmatrix} 
	    D        & \boldzero \\ 
	   \boldzero & \boldzero 
	\end{psmallmatrix}
    } \; = \; 
    (
	\boldz^{
	    \begin{psmallmatrix} 
		D \\ 
		\boldzero  
	    \end{psmallmatrix}
	},
	\boldz^{
	    \begin{psmallmatrix} 
		\boldzero \\ 
		\boldzero 
	    \end{psmallmatrix}
	}
    ) \; = \; 
    \boldb^Q \; = \; 
    (\boldb^{Q_r}, \boldb^{Q_0}) \;.
\]
Since $D = \op{diag} (d_1,\dots,d_r) \in M_{r \times r}(\Z)$,
the original system \,$\boldx^A = \boldb$\, can now be decomposed into a combined system
\begin{align}
    (z_1,\dots,z_r)^{
	\begin{psmallmatrix}
	    d_1 &        &     \\
		& \ddots &     \\
		&        & d_r
	\end{psmallmatrix}
    } &= \boldb^{Q_r} \label{equ:diagonal}\\
    \boldone \; &= \; \boldb^{Q_0} \label{equ:consistency} \\
    z_{r+1},\dots,z_{n} &:\;\; \text{free} \label{equ:free-vars}
\end{align}
where \eqref{equ:consistency} appears when $r < m$ 
with $\boldone = (1,\dots,1) \in (\C^\ast)^{m-r}$,
and \eqref{equ:free-vars} appears when $r < n$. 
The word ``free'' in \eqref{equ:free-vars} means the system 
imposes no constraints on the $n-r$ variables $z_{r+1},\dots,z_{n}$.

Focusing on the above decomposed system,
it is clear that if $r < m$, then the system is inconsistent unless $\boldone=\boldb^{Q_0}$.
If the system is consistent (namely, \eqref{equ:consistency} holds),
then the solutions to \eqref{equ:diagonal} are exactly
\begin{equation}
    \label{equ:z-solns}
    \begin{cases}
        z_1 = e^{2 k_1 \pi / d_1} \zeta_1 &\text{for} \quad k_1 = 0,\dots,d_1 - 1 \\
	    z_2 = e^{2 k_2 \pi / d_2} \zeta_2 &\text{for} \quad k_2 = 0,\dots,d_2 - 1 \\
	    \vdots &\vdots \\
	    z_r = e^{2 k_r \pi / d_r} \zeta_r &\text{for} \quad k_r = 0,\dots,d_r - 1
    \end{cases}
\end{equation}
where each \,$\zeta_j$\, is a fixed choice of the \,$d_j$-th\, root of \,$j$-th\, coordinate of \,$\boldb^Q$.
Clearly, all of them are isolated and the total number of these solutions is $|\prod_{j=1}^r d_j| = |\det D|$.
If $r < n$, then the solution set of the decomposed system \eqref{equ:diagonal}--\eqref{equ:free-vars}
in $(\C^\ast)^n$ breaks into ``components'' of the form 
$\{ (e^{2k_1 \pi/d_1} \zeta_1,\dots,e^{2k_r \pi/d_r} \zeta_r,z_{r+1},\dots,z_{n}) \, : \, (z_{r+1},\dots,z_n) \in (\C^\ast)^{n-r} \}$,
and they are in one-to-one correspondence with solutions in \eqref{equ:z-solns}.
Since each component is parametrized by the $n-r$ free variables $z_{r+1},\dots,z_n$,
it is smooth and of dimension $n-r$. 
Furthermore, they are disjoint, because these components have distinct $z_1,\dots,z_r$ coordinates.

To translate the above description of the $(\C^\ast)^n$-solution set 
of the decomposed system (in $\boldz$) into a description the original
solution set $\V(\boldx^A - \boldb)$, 
one may simply apply the change of variables $\boldx = \boldz^P$.
Note that this map and its inverse $\boldz = \boldx^{P^{-1}}$ are both given by 
by monomials (\emph{bi-regular} maps \cite{hartshorne_algebraic_1977}),
the basic properties of the solution set, such as, the number of solution components,
their dimensions, and smoothness are therefore preserved.
To summarize, the above elaborations assert the following proposition.

\begin{proposition}[
    Global parametrization
    \cite{eisenbud_binomial_1996,kahle_decompositions_2010,sturmfels_equations_1997}] 
    For the solution set $\V(\boldx^A - \boldb)$ in $(\C^\ast)^n$,
    let $P,Q,Q_0$ and $D$ be those matrices appeared in the decompositions of $A$
    in \eqref{equ:smith} and \eqref{equ:rank-decomp},
    and let $r = \op{rank} A$.

    If \,$\boldone \ne \boldb^{Q_0}$\, then the binomial system is inconsistent,
    and hence its solution set in $(\C^\ast)^n$ is empty.

    If $\boldone = \boldb^{Q_0}$ then 
    the solution set of $\boldx^A = \boldb$ in $(\C^\ast)^n$ consists of
    $|\prod_{j=1}^r d_j| = |\det D|$ connected components 
    $V_{k_1,\dots,k_r}$ for $k_1 \in \{0,\dots,d_1-1\}, \dots, k_r \in \{0,\dots,d_r-1\}$.
    Each component $V_{k_1,\dots,k_r}$ is smooth of dimension $n-r$,
    and it is parametrized by the smooth global parametrization
    $\phi_{k_1,\dots,k_r} \,:\, (\C^\ast)^{(n-r)} \to V_{k_1,\dots,k_r}$ given by
    \begin{equation}
	\label{equ:parametrization}
	\phi_{k_1,\dots,k_r}(t_1,\dots,t_{n-r}) \; = \; (e^{2 k_1 \pi / d_1} \zeta_1,\dots,e^{2 k_r \pi / d_r} \zeta_r,t_1,\dots,t_{n-r})^P
    \end{equation}
    where each $\zeta_j$ is a fixed choice of the $d_j$-th root of the $j$-th coordinate of  \,$\boldb^Q$.
\label{pro:parametrization}
\end{proposition}

Note that, as previously stated, in the case of $r=n$,
the solution set $\V(\boldx^A-\boldb)$ is of dimension $n-r = 0$,
that is, $\V(\boldx^A-\boldb)$ consists of isolated points.
Then the ``parametrizations'' $\phi_{k_1,\dots,k_r}$ are understood as 
constants each describes a single isolated point.

As indicated in Proposition \ref{pro:parametrization}, 
for a consistent Laurent binomial system $\boldx^A = \boldb$ where $A \in M_{n \times m}(\Z)$ 
with $r = \op{rank}(A) < n$, each component of the solution set in $(\C^\ast)^n$ 
will be of dimension $n - r > 0$.
In this situation, for both theoretical interests and demands from concrete applications,
like the Master Space problem to be discussed in \S \ref{sec:mspace},
one often wishes to identify another important property: the \emph{degrees} of the components.
Degree is a classic concept developed for plane algebraic curves.
For example, the quadratic equation $y-x^2=0$\, defines a curve of degree 2, i.e., the parabola.
The generalized notation of degree for irreducible algebraic sets
is usually formulated algebraically via \name{Hilbert} Polynomials. 
In this article, following the common practice of Numerical Algebraic Geometry, 
we shall take a geometric approach:
Let $V = V_{k_1,\ldots,k_r}$ be a component of $\V(\boldx^A-\boldb)$
for some fixed choice of $k_1,\ldots,k_n$ as defined in Proposition \ref{pro:parametrization}.
The number of isolated intersection point between $V$ and
a ``generic'' affine space of complementary dimension is a fixed number,
and this number is the \term{degree} of $V$, denoted by $\deg V$.
In algebraic terms, we are considering the degree of the \emph{projective closure} of $V$.

Stated more precisely,
let $\mathcal{G}_r$ be the set of all affine space in $\C^n$ of dimension $r = n - \dim V$.
Then it can be shown that in a fixed open and dense subset of $\mathcal{G}_r$,
all the affine spaces intersect with $V$ at a fixed number of isolated points.
This geometric interpretation of degree is explained in 
\cite{fulton_intersection_1998,hartshorne_algebraic_1977,sommese_numerical_2005}.

From a computational standpoint,
a generic affine space in $\mathcal{G}_r$ can be represented by the solution set
of a system of $d := n - r$ linear equations with generic coefficients.
Therefore $\deg V$ is precisely the number of points 
$\boldx = (x_1,\ldots,x_n) \in V$ that satisfies the system of linear equations
\begin{equation} \label{equ:linear-cut}
    \begin{cases}
        c_{11} x_1 + c_{12} x_2 + \cdots + c_{1n} x_n &= c_{10} \\
                                                      &\mathrel{\makebox[\widthof{=}]{\vdots}}   \\
        c_{d1} x_1 + c_{d2} x_2 + \cdots + c_{dn} x_n &= c_{d0}.
    \end{cases}
\end{equation}
where $c_{ij}$ for $i=1,\ldots,d$ and $j=1,\ldots,n$
are generic complex numbers.
But recall that the set $V = V_{k_1,\ldots,k_r}$ is precisely the image of the injective map
\[ \phi_{k_1,\dots,k_r}(t_1,\dots,t_{d}) \; = \; (e^{2 k_1 \pi / d_1} \zeta_1,\dots,e^{2 k_r \pi / d_r} \zeta_r,t_1,\dots,t_{d})^P \]
in Proposition \ref{pro:parametrization}.
If we let $\xi = (e^{2 k_1 \pi / d_1} \zeta_1,\dots,e^{2 k_r \pi / d_r}$
and $\boldt = (t_1,\ldots,t_d)$ then
\[ 
    \phi_{k_1,\dots,k_r} (\boldt) = 
    (\xi, \boldt)^{ \left( \begin{smallmatrix} P_r \\ P_0 \end{smallmatrix} \right) } = 
    (\, \xi^{\boldp_r^{(1)}} \boldt^{\boldp_0^{(1)}}, \, \ldots, \, \xi^{\boldp_r^{(n)}} \boldt^{\boldp_0^{(n)}} \,)
\]
where for each $j = 1, \ldots, n$, $\boldp_r^{(j)}$ and $\boldp_0^{(j)}$
are the $j$-th columns of $P_r$ and $P_0$ respectively.
In other words, $V$ has the global parametrization $x_i = \xi^{\boldp_r^{(1)}} \boldt^{\boldp_0^{(1)}}$.
Therefore the intersections between $V$ and the generic affine space
defined by \eqref{equ:linear-cut} are precisely the solutions of the polynomial system
\begin{equation*}
    \begin{cases}
    	c_{11} \, \xi^{\boldp_r^{(1)}} \, \boldt^{\boldp_0^{(1)}} +
    	c_{12} \, \xi^{\boldp_r^{(2)}} \, \boldt^{\boldp_0^{(2)}} + \cdots + 
        c_{1n} \, \xi^{\boldp_r^{(n)}} \, \boldt^{\boldp_0^{(n)}} & =  c_{10} \\
        
%
    	                                                          & \mathrel{\makebox[\widthof{=}]{\vdots}} \\
                                                      
    	c_{d1} \, \xi^{\boldp_r^{(1)}} \, \boldt^{\boldp_0^{(1)}} +
        c_{d2} \, \xi^{\boldp_r^{(2)}} \, \boldt^{\boldp_0^{(2)}} + \cdots + 
        c_{dn} \, \xi^{\boldp_r^{(n)}} \, \boldt^{\boldp_0^{(n)}} & =  c_{d0}
    \end{cases}
\end{equation*}

By letting $c_{ij}' := c_{ij} \xi^{\boldp_r^{(j)}} \in \C$ and $c_{i0}' = c_{i0}$
for each $i=1,\dots,d$ and $j=1,\dots,n$,
the above is a systems of $d$ polynomial equations in variables 
$\boldt = (t_1,\dots,t_d)$ with generic complex coefficients $c_{ij}'$ 
and the same set of monomials $\boldt^{\boldp_0^{(j)}}$.

%

\begin{proposition}[Degree via affine space cut] \label{pro:cut}
    If $r < n$ and $\V(\boldx^A-\boldb) \ne \varnothing$,
    then the degree of each component $V$ of $\V(\boldx^A-\boldb)$
    agrees with the number of solutions $\boldt \in (\Cstar)^d$
    of the system of $d$ Laurent polynomial equations
    \begin{equation} \label{equ:affine-cut}
        \begin{cases}
            c_{11} \boldt^{\boldp_0^{(1)}} +
            c_{12} \boldt^{\boldp_0^{(2)}} + \cdots + 
            c_{1n} \boldt^{\boldp_0^{(n)}} & =  c_{10} \\
            & \mathrel{\makebox[\widthof{=}]{\vdots}} \\
            
            c_{d1} \boldt^{\boldp_0^{(1)}} +
            c_{d2} \boldt^{\boldp_0^{(2)}} + \cdots + 
            c_{dn} \boldt^{\boldp_0^{(n)}} & =  c_{d0}
        \end{cases}
    \end{equation}
    for generic complex coefficients $c_{ij} \in \C$.
\end{proposition}

It is important to note that for generic coefficients, the $\Cstar$-solutions of
the above system are all isolated (0-dimensional), and the total number is a constant.
Indeed, in \cite{kushnirenko_newton_1975}, \name{Kushnirenko} has shown that
this number can be expressed in terms of the volume of a geometric object:
the \emph{Newton polytope} of the above system.
Here we state the result in the context of degree computation
and leave the technical statement of \name{Kushnirenko}'s theorem,
as well as its related concepts to Appendix \S \ref{sec:kushnirenko}.

\begin{proposition}[Degree as volume] \label{pro:degree}
    For generic choices of the coefficients,
    the number of solutions $\boldt \in (\Cstar)^d$ of the
    system of Laurent polynomial equation \eqref{equ:affine-cut}
    and hence the degree of $V$ is
    \begin{equation} \label{equ:deg-vol}
        \deg V = d\, ! \cdot \op{Vol}_{d} (\op{conv} \{ \boldp_0^{(1)}, \dots, \boldp_0^{(n)}, \boldzero \})
    \end{equation}
    where  $\boldzero = (0,\dots,0)^\top \in \R^d$ 
    and columns \,$\boldp_0^{(1)},\dots,\boldp_0^{(n)}$ of the matrix $P_0$ 
    are considered as points in $\R^{d}$. 
    The notation $\op{conv}$ denotes the operation of taking convex hull,
    and $\op{Vol}_{d}$ is the volume of a convex body in $\R^d$.
\end{proposition}

The degree of the solution set of $\boldx^A = \boldb$ can be computed
efficiently through methods in combinatorial geometry. 

\section{Parallel Smith Normal Form computation} \label{sec:computation}

As summarized in \autoref{pro:parametrization},
the key to finding the dimension, number of components, and global parametrization
of the $\Cstar$-solution set $\V(\boldx^A-\boldb) \subset (\Cstar)^n$
is the Smith Normal Form \eqref{equ:smith} of the exponent matrix $A$.
In this section, we briefly review the procedure for computing 
the Smith Normal Form of an integer matrix and then outline the
parallel modification that is suitable for both multi-core systems and GPU.

We first briefly review the standard algorithm for computing 
the Smith Normal Form of a matrix with integer entries.
One of the classic algorithms for computing the Smith Normal Form uses
successive row ($n$) and column ($m$) reductions of the input matrix,
as listed in \cite[Algorithm 2.4.14]{cohen_course_1993}
and \cite[Section 8.5.1]{gockenbach_finite-dimensional_2011}:
Consider the special case where $A = \left( \begin{smallmatrix} a_1 \\ a_2 \end{smallmatrix} \right)$
with $a_1, a_2$ nonzero, that is, take $n = 2$ and $m = 1$. 
By the B\'ezout's identity, there exist $s$ and $t$ such that
$d := \op{gcd} (a_1, a_2) = s \, a_1 + t \, a_2$. Let
\[
    P =
    \begin{pmatrix}
        s              & t   \\
        -\frac{a_2}{d} & \frac{a_1}{d}
    \end{pmatrix} ,
\]
then $\det P = \frac{s a_1 + t a_2}{d} = \frac{d}{d} = 1$ and 
\[
    P A = 
    \begin{pmatrix}
    s              & t   \\
    -\frac{a_2}{d} & \frac{a_1}{d}
    \end{pmatrix}
    \begin{pmatrix}
        a_1 \\ a_2
    \end{pmatrix}
    =
    \begin{pmatrix}
        s a_1 + t a_2 \\
        - \frac{a_2 a_1}{d} + \frac{a_2 a_2}{d}
    \end{pmatrix}
    =
    \begin{pmatrix}
        d \\ 0
    \end{pmatrix}
\]
Similarly, for the special case $A = \left( \begin{smallmatrix} a_1 & a_2 \end{smallmatrix} \right)$,
let $Q = \left( \begin{smallmatrix} s & \frac{-a_2}{d} \\ t & \frac{a_1}{d} \end{smallmatrix} \right)$,
then $A \, Q = \begin{pmatrix} d & 0 \end{pmatrix}$.
In general, $n \times n$ and $m \times m$ version of of the above matrices $P$ and $Q$
can be constructed to perform row and column reduction respectively
for a $n \times m$ integer matrix.

After repeated such row and column reduction together with potential row and column permutations
one can construct unimodular matrices $P^{(1)}, \ldots, P^{(k)} \in \zmat{n}{n}$
and $Q^{(1)}, \ldots, Q^{(\ell)} \in \zmat{m}{m}$ such that
\[
    P^{(k)} \cdots P^{(1)} \, A \; Q^{(1)} \cdots Q^{(\ell)} = 
    \left(
    \begin{smallmatrix}
    	d_1 &        &     &   &        &  \\
    	    & \ddots &     &   &        &  \\
    	    &        & d_r &   &        &  \\
    	    &        &     & 0 &        &  \\
    	    &        &     &   & \ddots &  \\
    	    &        &     &   &        & 0
    \end{smallmatrix}
    \right)
\]
with $r = \op{rank} A$ and $d_1, \dots, d_r$ nonzero.
As noted in standard references such as \cite{gockenbach_finite-dimensional_2011},
further reduction can ensure $d_1 \mid d_2 \mid \cdots \mid d_r$, 
but for the purpose of solving binomial system it is not necessary.

By their design, GPUs are naturally well suited to perform the row and column reductions 
\cite{nvidia_corporation_nvidia_2011} used in computing the Smith Normal Form.
As \autoref{tab:gpu-row-reduction} shows, 
GPUs have a clear advantage over CPUs in performing simple row reductions
for sufficiently large matrices: over 30 fold speedup can be achieved.
\S \ref{sec:mspace} shows the result of this algorithm when applied to the master space problem.

\begin{table}
    \centering
    \begin{tabular}{|c|l|}
        \hline
        Row length   & Speedup ratio \\ \hline
        $  10$ &  $0.00$ \\ \hline
        $  50$ &  $0.00$ \\ \hline
        $ 100$ &  $1.91$ \\ \hline
        $ 200$ &  $1.99$ \\ \hline
        $ 400$ &  $8.01$ \\ \hline
        $ 800$ & $14.20$ \\ \hline
        $1600$ & $22.00$ \\ \hline
        $3200$ & $31.79$ \\ \hline
    \end{tabular}
    \vskip 1ex
    \caption{
        Speedup ratio achieved by a \techy{NVidia GTX 780} graphics card (GPU) 
        on a \emph{single} row reduction operation on when compared to an equivalent
        single-threaded CPU-based implementation on a \techy{Intel Core i7 4770k} CPU.
        In each case, the speedup ratio is computed as the average of three runs.
        Time consumed by data transfer is \emph{not} computed.
        }
    \label{tab:gpu-row-reduction}
\end{table}
\section{Parallel degree computation} \label{sec:degree}

When the solution set consists of positive dimensional components,
Proposition \ref{pro:degree} provides a computationally viable means for computing
the degree of each component as the normalized volume of a convex polytope.
In this section we shall present a parallel algorithm for computing the degree
that is suitable for both multi-core systems and GPUs, though the focus is the
GPU-based implementation.
Throughout this section, let $V$ be a component of 
$\V(\boldx^A-\boldb) \subset (\Cstar)^n$,
and let $d = \dim V = n - r$.
Here we shall focus on the case where $d > 0$.
Let $P^0 = ( \boldp_0^{(1)}, \ldots, \boldp_0^{(n)} ) \in \zmat{d}{n}$
be the matrix appears in the Smith Normal For of $A$ \eqref{equ:smith}.
Considering each $\boldp_0^{(j)}$ as a point in $\R^d$,
let $S = \{ \boldp_0^{(1)}, \ldots, \boldp_0^{(n)} \} \subset \R^d$ 
be the finite point set. Then by Proposition \ref{pro:degree},
\begin{equation} \label{equ:deg-vol2}
    \deg V = d \, ! \Vol_d (\conv S).
\end{equation}
For brevity, let $\nvol_d = d\,! \Vol$ be the \term{normalized volume}
function in $\R^d$, then the above equation can be written as
\begin{equation} \label{equ:deg-nvol}
    \deg V = \nvol_d (\conv S)
\end{equation}
Therefore the task of computing the degree for $V$ is equivalent to
the computation of the normalized volume of a \emph{lattice polytope}
(a polytope whose vertices have integer coordinates).

\begin{remark}
    Clearly, \eqref{equ:deg-vol2} and \eqref{equ:deg-nvol} are equivalent.
    However, from the computational point of view, there is one crucial distinction:
    The knowledge that $\nvol_d(S)$ must be an integer permits the use of efficient 
    but potentially less accurate numerical methods using floating point arithmetic 
    and still obtain the correct result.
    Indeed, the exact results can still be obtained as long as the 
    total absolute error is kept below $1/2$.
    This is not possible for methods that are designed to compute volume of
    more general polytopes.
    The algorithm, for computing \eqref{equ:deg-nvol}, to be presented below, 
    is hence not directly comparable to \emph{exact volume computation} algorithms
    \cite{barany_computing_1987,bueler_exact_2000} for general polytopes.
\end{remark}

Our parallel algorithm for the degree computation is developed based on the 
parallel ``mixed cell enumeration'' algorithm presented in \cite{chen_mixed_2014}.
(See Remark \ref{rmk:history} below)
Among many different approaches for computing the normalized volume,
here we adopt a technique known as \emph{regular simplicial subdivision} \cite{loera_triangulations:_2010}.
This approach produces an important byproduct that will be used in the computation
of witness set, which will be the subject of \S \ref{sec:witness}.
In this approach, we are interested in computing the normalized volume $\nvol_d(\conv S)$ 
by dividing the lattice polytope $\conv S$ into a collection of smaller pieces for which
the volume computation is easy.
\begin{definition} \label{def:subdiv}
    A \term{cell} of $S$ is simply an affinely independent subset of $S$.
    A \term{simplicial subdivision} of $S$ is a collection $\calD$ of cells of $S$, 
    such that
    \begin{enumerate}
        \item For each $C \in \calD$, $\conv C$ is a $d$-simplex inside $\conv S$;
        \item For any distinct pair of simplices $C_1, C_2 \in \calD$,
            the intersection of $\conv C_1$ and $\conv C_2$, if nonempty, 
            is a common face of the two; and
        \item The union of convex hulls of all cells in $\calD$ is exactly $\conv S$.
    \end{enumerate}
\end{definition}
A simplicial subdivision plays an important role in computing $\nvol_d(\conv S)$:
the normalized volume of a $d$-simplex in $\R^d$ is easy to compute:
given a $d$-simplex $\Delta = \conv \{ \bolda_0, \dots, \bolda_d \} \subset \Z^d$,
\begin{equation} \label{equ:simplex-vol}
    \nvol_d (\Delta) = \det
    \begin{pmatrix}
        \bolda_0 & \dots & \bolda_k
    \end{pmatrix}.
\end{equation} 
So the volume of $\conv S$ can be computed easily 
as the sum of the volume of all simplices in $\calD$.

Note that the simplicial subdivision for a given polytope is, in general, 
not unique, and there are many different approaches for constructing them.
Here we focus on the approach of regular simplicial subdivision:
One can define a ``lifting function'' $\omega : S \to \R$
by assigning a real number to each point in $S$.
For each point $\bolda \in S$,
a new point $(\bolda,\omega(\bolda)) \in \R^n$ can be created
by using $\omega(\bolda)$ as an additional coordinate.
This procedure ``lifts'' points of $S$ into $\R^{d+1}$,
the space of one higher dimension. Let
\begin{equation}
    \hat{S} = \{ \hat{\bolda} = (\bolda, \omega(\bolda)) \mid \bolda \in S \}
\end{equation}
be the lifted version of $S$ via the lifting function $\omega$.
Figures \ref{fig:lifting-a} and \ref{fig:lifting-b} show examples
of this lifting procedure..
Let $\pi : \R^{d+1} \to \R^d$ be the projection that simply erases the last coordinate,
then $\pi(\hat{S}) = \cal{S}$.

Recall that for a face $\hat{F}$ of the lifted polytope $\conv \hat{S}$, 
its \emph{inner normal} is a vector $\hat{\boldalpha} \in \R^{d+1}$ 
such that the linear functional $\langle \bullet, \hat{\boldalpha} \rangle$ 
attains its minimum over $\conv \hat{S}$ on $\hat{F}$.
Moreover, a face $\hat{F}$ of $\conv \hat{S}$ is called a \term{lower face} 
with respect to the projection $\pi$ if its inner normal $\hat{\boldalpha}$ 
has positive last coordinate.
Without loss of generality, in this case, we may assume the last coordinate 
of $\hat{\boldalpha}$ to be 1, that is, $\hat{\boldalpha} = (\nlist{\alpha},1) \in \R^{d+1}$.
It can be shown that for \emph{almost all} choices of the lifting function $\omega : S \to \R$,
the projections of all the $d$-dimensional lower faces 
of $\conv \hat{S}$ via $\pi$ form a simplicial subdivision for $\conv S$
which is called a \emph{regular simplicial subdivision} of $\conv S$.
The construction of this simplicial subdivision is therefore equivalent to 
the enumeration of all the lower faces of $\conv \hat{S}$.

\begin{example}
    Consider, for example, $S = \{(0,0),(0,1),(1,1),(1,0)\}$ in the $xy$-plane.
    A simplicial subdivision of $\conv S$ can be obtained via the following procedure:
    First one assign ``liftings'' $\omega_1,\omega_2,\omega_3,\omega_4 \in \R$
    to each of the vertices as the $z$-coordinate and obtain new points
    $(0,0,\omega_1),(0,1,\omega_2),(1,1,\omega_3),(1,0,\omega_4)$ in $\R^3$.
    It is easy to verify that with almost all choices of the liftings the four 
    ``lifted'' points (Figure \ref{fig:lifting-a}) do not lie on the same plane.
    In that case, the convex hull
    $\conv\{(0,0,\omega_1),(0,1,\omega_2),(1,1,\omega_3),(1,0,\omega_4)\}$ 
    of the four lifted vertices form a three dimensional polytope
    (Figure \ref{fig:lifting-b}) with triangle faces.
    Of particular importance is the lower hull of this polytope
    which are the faces facing downward.
    As shown in Figure \ref{fig:lifting-c}, the projection of the faces in the lower hull
    back onto the $xy$-plane form a simplicial subdivision of the original shape $\conv S$.
    
    \begin{figure}[h]
        \centering
        \begin{subfigure}{.33\textwidth}
            \centering
            \includegraphics{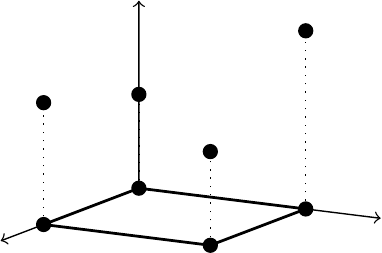}
            \caption{Lifted set $\hat{S}$}
            \label{fig:lifting-a}
        \end{subfigure}%
        \begin{subfigure}{.33\textwidth}
            \centering
            \includegraphics{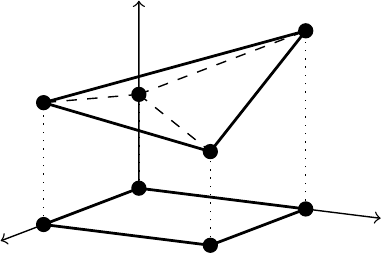}
            \caption{Lifted polytope $\conv \hat{S}$}
            \label{fig:lifting-b}
        \end{subfigure}%
        \begin{subfigure}{.33\textwidth}
            \centering
            \includegraphics{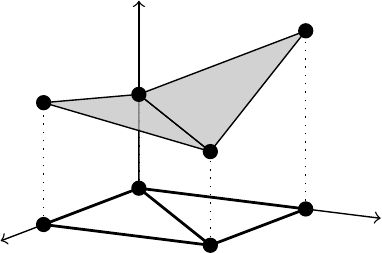}
            \caption{Projection of lower hull}
            \label{fig:lifting-c}
        \end{subfigure}
        \caption{Regular simplicial subdivision via generic lifting}
        \label{fig:lifting}
    \end{figure}
\end{example}

Algebraically speaking, a $d$-dimensional lower face of $\conv S$ is the convex hull
of a set of $d+1$ points $\{ \hat{\bolda}_0, \ldots, \hat{\bolda}_d \} \subset \hat{S}$
for which there exists a $\hat{\boldalpha} = (\boldalpha,1) \in \R^{d+1}$ such that
the system of inequalities
\begin{equation} \label{equ:lower-face}
    I(\bolda_0,\dots,\bolda_d) \; : \;
    \left\{
	\begin{aligned}
		\langle \hat{\bolda}_0, \hat{\boldalpha} \rangle &=
		\langle \hat{\bolda}_j, \hat{\boldalpha} \rangle & \text{ for } & j = 1, \dots, d \\
		\langle \hat{\bolda}_0, \hat{\boldalpha} \rangle &\le
		\langle \hat{\bolda}  , \hat{\boldalpha} \rangle & \text{ for } &\bolda \in S \;
	\end{aligned}
    \right.
\end{equation}
is satisfied.
In other words, the existence of the lower face defined by  
$\{ \hat{\bolda}_0, \ldots, \hat{\bolda}_d \} \subset \hat{S}$
is equivalent to the feasibility of the above system of inequalities $I(\bolda_0,\dots,\bolda_d)$.
This algebraic description of the lower faces is the basis on which enumeration methods are developed.
While a brute-force approach of checking all the possible combinations of $d+1$ points in $S$
against the system of inequalities \eqref{equ:lower-face} may be possible, 
the combinatorial explosion will likely render it impractical for all but the most trivial cases.

In the following subsections, we shall present an approach 
that results in a \emph{parallel algorithm} 
which is suitable for both multi-core systems and GPU devices.
In this approach, we employ two complementing processes of ``extension'' and ``pivoting''.
We shall outline them below.

\begin{remark}[Connection to existing works] 
    \label{rmk:history}
    The approach developed in this article is a natural continuation of a rich web of works
    on the ``mixed cell enumeration'' problem initiated by the seminal work \cite{huber_polyhedral_1995}.
    The degree computation problem can be considered as a special case of the mixed cell enumeration problem,
    and the connection is explained in \S \ref{sec:kushnirenko}.
    Active development in the algorithmic aspects of this problem can be found in works such as
    \cite{chen_mixed_2014,gao_mixed_2000,gao_mixed_2003,lee_mixed_2011,li_finding_2001,mizutani_demics:_2008,mizutani_dynamic_2007,verschelde_mixed-volume_1996}.
    A broad survey of this topic can be found in \cite{li_numerical_2003}.
    
    The ``pivoting'' process (for ``mixed cell enumeration''), 
    to be described below, was proposed in \cite{gao_mixed_2000}.
    However, in terms of performance, it was quickly eclipsed by the ``extension'' process 
    developed in \cite{li_finding_2001}, \cite{gao_mixed_2003}, and \cite{mizutani_demics:_2008,mizutani_dynamic_2007}.
    In the present article, the pivoting and the extension processes are combined
    as we believe the complementing duo offers much better scalability which
    is crucial in the GPU-based implementations.
    This is confirmed by the numerical experiments, to be presented in \S \ref{sec:mspace}.
    
    The graph-theoretic view of the ``cell enumeration'' process,
    adopted in this article, was originally developed  in \cite{gao_mixed_2000} 
    and, independently, in \cite{mizutani_dynamic_2007}.
    The parallelization of the algorithm follows the same general idea attempted 
    in \cite{chen_mixed_2014}, but it is modified, in this article, 
    to adapt to the massively parallel GPU architectures.
\end{remark}

\subsection{Extension of $k$-faces} \label{sec:extension}

Intuitively speaking, in the extension process, 
one starts with the vertices of the lower hull of $\conv \hat{S}$.
For each of these vertices, systematic attempts are made to ``extend'' it 
by finding another lower vertex so that the two vertices form a ``lower edge''
(an edge on the lower hull of $\conv \hat{S}$).
The possible extensions may not be unique, and for each possibility,
further attempts are made to extend it to 2-dimensional lower faces.
This process continues until one reaches all the $d$-dimensional lower faces.
Finally, the collection of such $d$-dimensional lower faces will project down,
via $\pi$, to form a simplicial subdivision for $\conv S$.

To describe this process, we first extend the characterization 
\eqref{equ:lower-face} to include lower faces of all dimensions:
A set of affinely independent $k+1$ points in $\hat{S}$ is said to determine a 
\term{lower $k$-face} if their convex hull form a $k$-dimensional lower face
of $\conv \hatS$ with respect to the projection $\pi$.
Stated algebraically, the affinely independent set $\{ \bolda_0,\dots,\bolda_k \}$
determines a lower $k$-face if and only if there exists 
an $\hat{\boldalpha} = (\boldalpha,1) \in \R^{d+1}$,
such that the system of inequalities
\begin{equation} \label{equ:lower-k-face}
    I(\bolda_0,\dots,\bolda_k) \; : \;
    \left\{
	\begin{aligned}
		\langle \hat{\bolda}_0, \hat{\boldalpha} \rangle &=
		\langle \hat{\bolda}_j, \hat{\boldalpha} \rangle & \text{ for } & j = 1, \dots, k \\
		\langle \hat{\bolda}_0, \hat{\boldalpha} \rangle &\le
		\langle \hat{\bolda}  , \hat{\boldalpha} \rangle & \text{ for } &\bolda \in S \;
	\end{aligned}
    \right.
\end{equation}
is satisfied.

\begin{figure}
    \centering
    \includegraphics{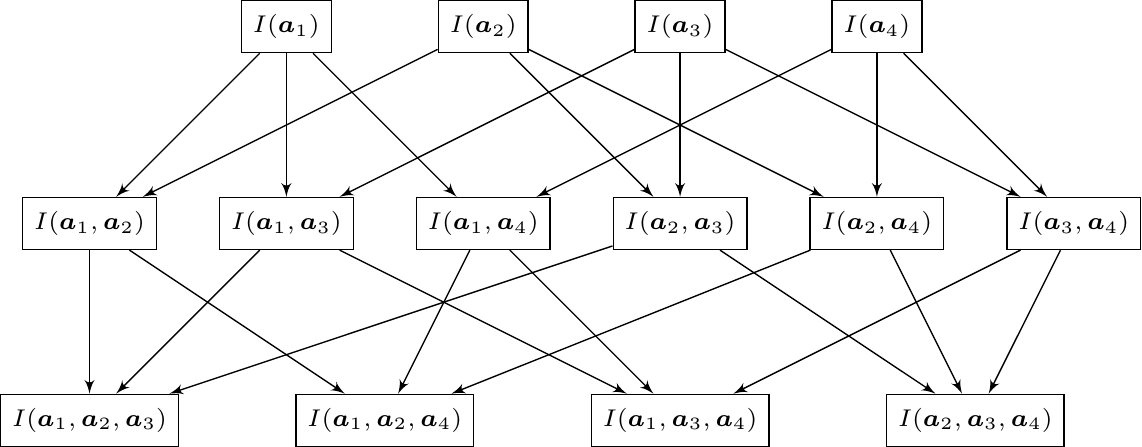}
    \caption{A directed acyclic graph of possible lower $k$-faces}
    \label{fig:dag}
\end{figure}

Clearly, a lower $0$-face is a vertex on the lower hull of $\conv \hat{S}$.
Similarly, a lower $1$-face is simply a lower edge.
We can conveniently organize all possible system of inequalities
of the above form into a \emph{directed acyclic graph}, 
as illustrated ine \autoref{fig:dag},
where each node represents a system of inequalities and
there is an edge from $I(A)$ to $I(B)$ whenever $B$ is obtained by joining 
a new points in $S$ into $A$.
With this construction, the resulting graph is \emph{graded}
by the number of points involved.

It can be easily verified that for generic lifting function $\omega$,
containment relation between lower $k$-faces of the same dimension is impossible.
That is, for a fixed $k$, no lower $k$-face is contained in another lower $k$-face.
Therefore the the graph describes precisely the containment relationship among
possible lower $k$-faces.

\begin{figure}
    \centering
    \includegraphics{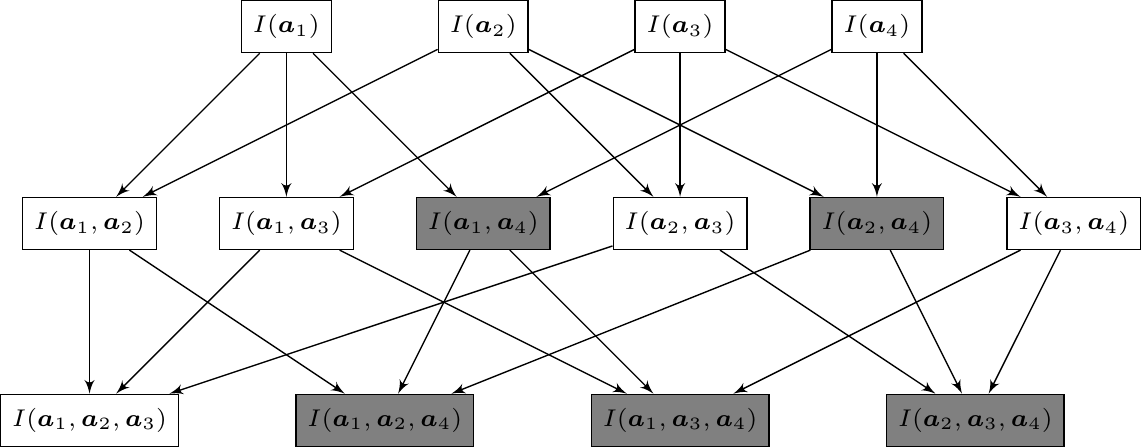}
    \caption{
        A direct graph of possible lower $k$-faces
        colored by the feasibility of the corresponding system of inequalities.
        }
    \label{fig:dag-colored}
\end{figure}

A node is said to be feasible if the corresponding system of inequality is feasible.
\autoref{fig:dag-colored} shows an example of the labeling of the graph
via the feasibility of the nodes: dark for infeasible nodes and white for feasible ones.
Recall that a node determines to lower $k$-faces if and only if it is feasible.
Hence we only need to explore of the feasible subgraph
(the white subgraph in \autoref{fig:dag-colored}).

One crucial observation is that if two points do not define a lower edge,
then they cannot be a part of any lower faces of dimension greater than 1.
More generally, if a set of points does not define a lower $k$-face,
then there are no lower $j$-faces containing them for any $j > k$.
Stated formally, for $\hat{F_1} \subset \hat{S}$,
\begin{equation}
    I(\hat{F}_1) \text{ is infeasible } \implies I(\hat{F}) \text{ is infeasible for all } 
    \hat{F}_1 \subset \hat{F} \subset \hat{S}.
\end{equation}
In terms of the graph, if a node is infeasible,
then the entire subgraph reachable by that node is infeasible.

Therefore during the exploration of the graph,
once an infeasible node is encountered,
no further exploration from that node is needed
as all nodes reachable are infeasible.
This simple observation produces significant savings in terms of computation.

A key procedure in the exploration of the feasible subgraph
is the jump from one feasible node to another along an edge.
Assuming, for some $\{ \bolda_0,\dots,\bolda_k \} \subset S$, 
the node $I(\bolda_0,\dots,\bolda_k)$ is feasible,
then the feasibility of an adjacent node, say via the edge $\bolda_{k+1}$,
can be determined by solving the \emph{linear programming problem}
\begin{equation} \label{equ:ext-lp}
    LP(\bolda_0,\dots,\bolda_k \, ; \bolda_{k+1}) \; : \;
	\begin{gathered}
	    \text{Minimize } 
	    \langle \hat{\bolda}_{k+1}, \hat{\boldalpha} \rangle -
	    \langle \hat{\bolda}_0, \hat{\boldalpha} \rangle 
	    \text{ subject to } \\ 
    	\begin{aligned}
    		\langle \hat{\bolda}_0, \hat{\boldalpha} \rangle &=
    		\langle \hat{\bolda}_j, \hat{\boldalpha} \rangle & \text{ for } & j = 1, \dots, k \\
    		\langle \hat{\bolda}_0, \hat{\boldalpha} \rangle &\le
    		\langle \hat{\bolda}  , \hat{\boldalpha} \rangle & \text{ for } & \text{all } \bolda \in S \;
    	\end{aligned}
	\end{gathered}
\end{equation}
with the variable $\hat{\boldalpha} = (\boldalpha,1)$ for $\boldalpha \in \R^d$.

Note that under the constraints, the value of the objective function must be nonnegative.
Indeed, the minimum value of 0 is attainable precise when there is an $\hat{\boldalpha}$
for which the constraints are satisfied, simultaneously to
$\langle \hat{\bolda}_{k+1}, \hat{\boldalpha} \rangle = \langle \hat{\bolda}_0, \hat{\boldalpha} \rangle$.
That is, minimum value is 0 if and only if $I(\bolda_0,\dots,\bolda_k,\bolda_{k+1})$ is feasible.
In this case, a new node $I(\bolda_0,\dots,\bolda_k,\bolda_{k+1})$ is discovered.
Geometrically, we have ``extended'' the lower $k$-face determined by 
$\{ \hat{\bolda}_0,\dots,\hat{\bolda}_k \}$ into a lower $(k+1)$-face
by joining it the new vertex $\hat{\bolda}_{k+1}$.

Using the extension procedure as a basic building block,
we shall discuss, in \S \ref{sec:traverse}, we shall discuss the complete
parallel algorithm for the exploration of the feasible subgraph.

\subsection{Simplicial pivoting} \label{sec:pivot}

In the above we have described a process that gradually explore 
the feasible subgraph via extension procedures.
This process is complemented by another process which we shall call ``simplicial pivoting''
which explore the feasible subgraph by ``moving sideways'' in the graph
from one lower $d$-face to another.

\begin{wrapfigure}{r}{0.35\textwidth}
    \begin{center}
        \includegraphics{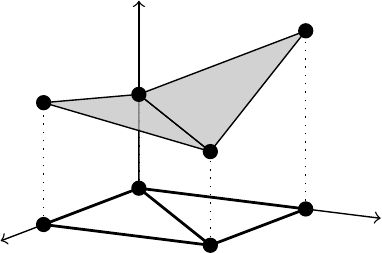}
    \end{center}
    \caption{Simplicial pivoting procedure moves from one lower face to another}
    \label{fig:flip}
\end{wrapfigure}

This process starts with a lower $d$-face of $\conv \hatS$ already obtained.
Consider, for example, one of the lower face shown in \autoref{fig:flip}.
Using an edge as a hinge, we shall ``pivot'' one lower face 
until another lower face is obtained.
More generally, recall that a lower $d$-face is determined by a set of $d+1$ 
affinely independent points $\{ \hat{\bolda}_0, \ldots, \hat{\bolda}_d \}$
in $\hatS$ that has an inner normal of the form $\hat{\boldalpha} = (\boldalpha,1)$
with $\boldalpha \in \R^{d+1}$. Stated algebraically, the system of inequalities
\begin{equation}
    \label{equ:lower-d-face}
    \begin{aligned}
        \langle \hat{\bolda}_0, \hat{\boldalpha} \rangle &=
        \langle \hat{\bolda}_j, \hat{\boldalpha} \rangle & \text{ for } & j = 1, \dots, d \\
        \langle \hat{\bolda}_0, \hat{\boldalpha} \rangle &\le
        \langle \hat{\bolda}  , \hat{\boldalpha} \rangle & \text{ for } & \text{all } \bolda \in S \;
    \end{aligned}
\end{equation} is satisfied.
Note that the $d$ equalities form a system of $d$ linearly independent constrains
on $\boldalpha \in \R^d$ and hence uniquely determines $\boldalpha$.
By removing a single equality from the above system,
we give the inner normal $\hat{\boldalpha}$ one degree of freedom
which would allow it to ``pivot''.
The goal is to let it pivot until it defines a different lower $d$-face.

For any choice $i = 0, \ldots, d$, with the equality corresponding to 
$\hat{\bolda}_i$ in the above system \eqref{equ:lower-d-face} removed,
the inner normal $\hat{\boldalpha} = (\boldalpha,1) \in \R^d$, 
now with one degree of freedom,  is characterized by the system
\begin{equation} \label{equ:pivot}
    P(\bolda_0,\dots,\bolda_d \, ; i) \; : \;
    \left\{
   	\begin{aligned}
   		\langle \hat{\bolda}_0, \hat{\boldalpha} \rangle &=
   		\langle \hat{\bolda}_j, \hat{\boldalpha} \rangle & \text{ for } & j = 1, \dots, d, \text{ but } j \ne i \\
   		\langle \hat{\bolda}_0, \hat{\boldalpha} \rangle &\lneq
   		\langle \hat{\bolda}_i, \hat{\boldalpha} \rangle \\
   		\langle \hat{\bolda}_0, \hat{\boldalpha} \rangle &\le
   		\langle \hat{\bolda}  , \hat{\boldalpha} \rangle & \text{ for } & \text{all } \bolda \in S \;
   	\end{aligned}
    \right.
\end{equation}
Note that this system has $d-1$ equalities.
If a solution with $d$ equalities exists,
then that solution corresponds to a different lower $d$-face.
In the context of Linear Programming, such a solution is called
a \emph{basic feasible solution}.
The problem of finding a basic feasible solution is known
as the \emph{Phase One} problem in Linear Programming.
It can be solved exactly and efficiently.

\begin{figure}[h]
    \centering
    \begin{subfigure}{0.32\textwidth}
        \centering
        \includegraphics{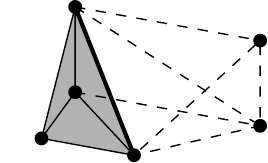}
        \caption{A simplex}
    \end{subfigure} %
    ~
    \begin{subfigure}{0.32\textwidth}
        \centering
        \includegraphics{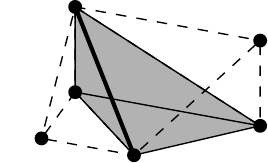}
        \caption{A simplex obtained by a single pivoting operation}
    \end{subfigure}%
    ~
    \begin{subfigure}{0.32\textwidth}
        \centering
        \includegraphics{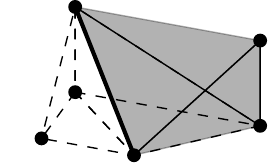}
        \caption{Another simplex obtained through further pivoting operation}
    \end{subfigure}
    \caption{
        Via the pivoting procedure, one moves from a simplex (a)
        to a different simplex (b) by leaving a chosen vertex
        and then to yet another simplex (c).
        }
    \label{fig:flip-3d}
\end{figure}

This procedure is called ``simplicial pivoting''.
It allows us to pivot from one lower $d$-face to another.
By repeatedly applying this procedure,
more lower $d$-faces can be gathered.
\autoref{fig:flip-3d} illustrates this process.

\subsection{Traverse the feasible subgraph} \label{sec:traverse}

In the above we have formulated the enumeration of lower $d$-faces
as the problem of exploring the feasible subgraph of which the
lower $d$-faces is a subset.
We also have two procedures for ``walking'' within the graph:
The extension procedure moves from one lower face to another of one higher dimension
while the pivoting procedure jumps from one lower $d$-face to another lower $d$-face.
With these building blocks in place,
the exploration can be handled by classic graph traversal algorithms
which we shall briefly review for completeness.

Most graph traversal algorithms follow a ``discover-explore''
procedure with proper book keeping \cite{skiena_algorithm_2009}:
They gradually explore the graph node by node through the connection 
between them while keeping track of the nodes visited so that no node
is explored twice.
For a single node, such an algorithm is divided into the \emph{discover} and \emph{explore} stages:
a node is first discovered, and then its connections to other yet unknown nodes are explored.
Clearly, each node only needs to be visited once.
That is, one only needs to explore a \emph{spanning tree} of the graph,
(a subgraph that contains all the vertices but is a tree in structure),
so some mechanism must be used to prevent a node from being visited twice.
To keep track of the nodes as they are being visited,
we assign each task a dynamic marker\,-- its state.
A node can be in one of the following three states:

\begin{description}
    \item[undiscovered] The initial status of every node.
    In this state, the existence of the node is completely unknown to us.
    
    \item[discovered] The existence of the node is known,
    but its connections to other nodes are not yet explored.
    
    \item[completely-explored] The existence of the node is known
    and its connections to other nodes have been fully explored.
\end{description}

Obviously, a node cannot be \emph{completely-explored} before it is first
discovered, so in the course of the algorithm, the state of vertices progresses
from \emph{undiscovered} to \emph{discovered} to \emph{completely-explored}.
This point of view also reveals the parallelism in such algorithms:
nodes on different branches of the spanning tree can be explored in parallel,
while consecutive nodes on a single branch must be discovered and explored in order.
To start the algorithm, an initial set of nodes are generated by some other means (bootstrapping).
The algorithm then discovers other nodes through their connections.
From these newly discovered nodes the algorithm can discover even more node.
This will continue as a self-sustaining process until all connected vertices are visited.

A complete algorithm also need a data structure to keep track of the discovered
but not yet completely explored vertices (bookkeeping).
In the present work, a \emph{queue} is used.
It is a linear data structure where newly discovered nodes are added 
to the back-end of the queue.
The use of the queue structure essentially imposes an implicit ordering of
``first-in-first-out'', that is, nodes discovered first are explored first.
In the context of graph traversal algorithms, this is referred to as
a \emph{breadth-first} strategy in exploring the feasible subgraph.
Experiments, presented in \cite{chen_mixed_2014}, suggest that
a more flexible ordering of nodes within the queue may provide better
performance, scalability, and memory usage.
However, for simplicity, in this work, only the breadth-first approach
has been studied.
The detail of this class of algorithms can be found in standard textbooks such as \cite{skiena_algorithm_2009}.
In \S \ref{sec:algo}, we list the pseudo code.

\subsection{Checking for duplicated discovery}

One important problem we must deal with, in the parallel algorithm, is that
same nodes may be discovered by different threads at the same time.
Since the degree is the sum of the normalized volume of the projection of all
the lower $d$-faces which are represented by the nodes in the feasible subgraph,
duplicated nodes will produce incorrect results.
Therefore, the mechanism for ensuring no duplicated lower $d$-faces are listed
is the key to the correctness of the algorithm.

This mechanism appears to be the bottleneck, in terms of performance, 
of the original algorithm \cite{gao_mixed_2000} for enumerating ``mixed cells'' 
using mainly the pivoting process which is one of the main inspiration of 
the present work (see Remark \ref{rmk:history}).
Our experiments confirm that an inefficient checking mechanism would be
the limiting factor of the scalability in a parallel implementation.
Since on a GPU, it is typical to have thousands of threads active simultaneously,
the efficiency of such a mechanism is crucial.

In the present work, the \emph{hash table} data structure is used to keep track of 
the nodes, in the graph, that have been discovered or completely-explored.
The great advantage of this choice is that unlike a sorted data structure,
hash table provides nearly constant access time, in \emph{most} cases.
In our current implementation, for simplicity, the well-tested 
bit-string hash function from the standard \techy{C++} library is used.

Our experiments suggest that a hash table with $2^{16}$ to $2^{20}$ entries
is sufficient for all problems considered in \S \ref{sec:mspace} in the sense 
that the collision rate in hash table access can be virtually ignored.


\subsection{Summary of the algorithm} \label{sec:algo}

In the above, we formulate the degree computation for solution components
defined by binomial systems as the exploration of the feasible subgraph
to be accomplished by the two complementing processes: extension and pivoting.
In this section, we list the main algorithms.

These algorithms are designed for a system with one or more GPU devices
and a single CPU with the GPU performance most of the computation intensive tasks.
For simplicity, we restrict ourselves to modern GPUs manufactured by \techy{NVidia}
and build our program based on \techy{NVidia CUDA} (a GPU programming framework).
All the GPU devices must share memory since they must all have access to
data structures \texttt{WaitingNodes}, \texttt{KnownNodes}, and \texttt{NewNodes}.
In the current implementation, this is accomplished via a technique 
known as \emph{pinned memory} \cite{nvidia_corporation_nvidia_2011} provided by the \techy{CUDA} framework.

In the following algorithms, the list \texttt{WaitingNodes} contains nodes 
whose feasibility are to be determined by the extension procedure.
\texttt{Cells} is the unordered collection of lower $d$-faces already discovered.
\texttt{KnownNodes} is the hash table that record the discovery of nodes,
and it is crucial mechanism by which we ensure the uniqueness of the discovered nodes.
Finally, \texttt{NewNodes} is an unordered list that keeps track of nodes discovered
through pivoting or extension.
They need to be checked against \texttt{KnownNodes} for uniqueness.

\textsc{Random} is a function that randomly choose an item from a collection
using pseudo random number generator.
The randomness is employed to achieve a more uniform performance 
from one run to another which simplifies the benchmarking process.
\textsc{SimplexPhaseOne} and \textsc{SimplexPhaseTwo} are the 
phase-one and phase-two algorithm of the simplex method for the linear programming 
problems \eqref{equ:ext-lp} and \eqref{equ:pivot} respectively.
Even though, at over 3000 lines, the \techy{C++} code for these two components
are the longest and most complicated parts of the entire program,
they have been a fixture of the long line of ``mixed volume computation''
software developed over the last two decades whence the present work inherits
much of its techniques and design.
Therefore we choose to not describe them in detail
and refer to works including \cite{chen_mixed_2014,gao_mixed_2000,gao_mixed_2003,lee_mixed_2011,li_finding_2001,mizutani_demics:_2008,mizutani_dynamic_2007}.

The \textsc{Extend} procedure tests the feasibility of a node (see \S \ref{sec:extension})
in the waiting list \texttt{WaitingNodes}, 
and it is designed to run simultaneously on all available threads across all GPU devices.

\begin{algorithmic}[1]
    \Function{Extend}{}
        \If {\texttt{WaitingQueue} $\ne \varnothing$}
            \State $\{ \bolda_0, \ldots, \bolda_k \} \gets $ \Call{Dequeue}{ \texttt{WaitingQueue} }
            \State $F \gets $ \Call{SimplexPhaseTwo}{ $LP(\{ \bolda_0, \ldots, \bolda_k \})$ }
            \If {$F \ne \varnothing$}
                \State \texttt{NewNodes} $\gets \texttt{NewNodes} \cup \{ F \}$ 
            \EndIf
        \EndIf
    \EndFunction
\end{algorithmic}

The \textsc{Pivot} procedure implements the simplicial pivoting process 
detailed in \S \ref{sec:pivot}, and it is designed to run simultaneously
on all available threads across all GPU devices.
It picks a random lower $d$-face already discovered 
and apply simplicial pivoting to potentially obtain a new lower faces.
Just like the \textsc{Extend} procedure above, 
newly discovered nodes will be place in the \texttt{NewNodes} list.

\begin{algorithmic}[1]
    \Function{Pivot}{}
        \If {\texttt{Cells} $\ne \varnothing$}
            \State $\{ \bolda_0, \ldots, \bolda_d \} \gets $ \Call{Random}{ \texttt{Cells} }
            \State $\ell \gets \min(d+1,10)$
            \For{$i = 1, \ldots, \ell$}
                \State $j \gets $ \Call{Random}{ $\{0,\ldots,d\}$ }
                \State $F \gets $ \Call{SimplexPhaseOne}{ $P(\{ \bolda_0, \ldots, \bolda_d \} \setminus \{ \bolda_j \})$ }
                \If {$F \ne \varnothing$}
                    \State \texttt{NewNodes} $\gets \texttt{NewNodes} \cup \{ F \}$ 
                \EndIf
            \EndFor
        \EndIf
    \EndFunction
\end{algorithmic}

The procedure \textsc{CheckUniq} checks newly discovered nodes against the hash table 
\texttt{KnownNodes} to make sure they have not already been discovered.
It  will run on a GPU device with a large number of threads simultaneously
checking the uniqueness of all nodes in the list of \texttt{NewNodes}.

\begin{algorithmic}[1]
    \Function{CheckUniq}{}
        \If {\texttt{NewNodes} $\ne \varnothing$}
            \State $\{ \bolda_0, \ldots, \bolda_k \} \gets $ \Call{Dequeue}{ \texttt{NewNodes} }
            \If { $\{ \bolda_0, \ldots, \bolda_k \} \not\in \texttt{KnownNodes}$ }
                \State $\texttt{KnownNodes} = \texttt{KnownNodes} \cup \{ F \}$
                \If { $k = d+1$ }
                    \State $\texttt{Cells} = \texttt{Cells} \cup \{ \{ \bolda_0, \ldots, \bolda_k \} \}$
                \Else
                    \ForAll { $\bolda \in S \setminus \{ \bolda_0, \ldots, \bolda_k \}$ }
                        \If { $\{ \bolda_0, \ldots, \bolda_k, \bolda \} \not\in \texttt{KnownNodes}$ }
                            \State $\texttt{WaitingQueue} = \texttt{WaitingQueue} \cup \{ \{ \bolda_0, \ldots, \bolda_k , \bolda \} \}$
                        \EndIf
                    \EndFor
                \EndIf
            \EndIf
        \EndIf
    \EndFunction
\end{algorithmic}

Finally, the main procedure, which runs on the CPU, 
coordinates all the different processes.
\begin{algorithmic}[1]
    \Function{Main}{}
        \State \texttt{WaitingNodes} $\gets S$
        \While { \texttt{WaitingNodes} $\neq \varnothing$ }
            \State Run \textsc{Extend} on available GPU threads
            \State Run \textsc{Pivot} on available GPU threads
            \State Wait for \textsc{Extend} and \textsc{Pivot}
            \State Run \textsc{CheckUniq} on available GPU threads
        \EndWhile
    \EndFunction
\end{algorithmic}

\section{Computation of witness sets} \label{sec:witness}

The concept of ``witness sets'' \cite{sommese2001numerical,sommese_numerical_2005}
is one of the most fundamental and versatile tool in numerical algebraic geometry.
In its most basic form, given a pure dimensional algebraic set,
it can be shown that its intersection with a generic affine space
of complementary dimension consists of finitely many isolated points.
This finite set is called a \emph{witness set} of the algebraic set.
It can be used to compute, among many other things, 
the irreducible decomposition and primary decomposition numerically.
In many scenarios, it produces the degree of each component as a byproduct.
Indeed, this technique (via witness sets) was first used to numerically compute 
the degrees of the ``Master Space'' problem in the work \cite{Hauenstein:2012xs}.

Given the ubiquity of the use of witness sets in numerical algebraic geometry,
in this section, we shall briefly outline a homotopy construction for
computing witness sets for a component of $\V(\boldx^A-\boldb)$.
It is a special case of the \emph{polyhedral homotopy} \cite{huber_polyhedral_1995}.

Recall that by Proposition \ref{pro:cut}, the intersection between a component
$V \subseteq \V(\boldx^A-\boldb)$ and a generic affine space of complementary
dimension consists of precisely the points $\boldt = (t_1,\ldots,t_d) \in (\Cstar)^d$
that satisfy the system of $d$ Laurent polynomial equation in $d$ variables given by
\begin{equation} \label{equ:affine-cut-again}
    \begin{aligned}
        c_{11} \boldt^{\boldp_0^{(1)}} +
        c_{12} \boldt^{\boldp_0^{(2)}} + \cdots + 
        c_{1n} \boldt^{\boldp_0^{(n)}} &=  c_{10} \\
        & \mathrel{\makebox[\widthof{=}]{\vdots}} \\
        c_{d1} \boldt^{\boldp_0^{(1)}} +
        c_{d2} \boldt^{\boldp_0^{(2)}} + \cdots + 
        c_{dn} \boldt^{\boldp_0^{(n)}} &=  c_{d0}
    \end{aligned}
\end{equation}
where the coefficients depends on both the choice of the component in 
$\V(\boldx^A-\boldb)$ and the choice of the $r$-dimensional affine space.

Reusing the notations from \S \ref{sec:degree},
let $S = \{ \boldp_0^{(1)}, \ldots, \boldp_0^{(n)} \}$,
and let $\omega : S \to \R$ be the generic lifting function used for constructing
regular simplicial subdivision of $\conv S$ in \S \ref{sec:degree}.
Without loss of generality, we can pick $\omega$ to have
images only in $\mathbb{Q}$.
With these, we introduce a new variable $s$ and consider
\begin{equation} \label{equ:affine-cut-homotopy}
    \small
    H (\boldt,s) = 
    \left\{
    \begin{matrix}
        c_{11} \boldt^{\boldp_0^{(1)}} s^{\omega(\boldp_0^{(1)})} + \cdots + 
        c_{1n} \boldt^{\boldp_0^{(n)}} s^{\omega(\boldp_0^{(n)})} -  
        c_{10} s^{\omega(\boldp_0^{(1)})} 
        = \sum_{\bolda \in S} c_{1,\bolda} \boldt^{\bolda} s^{\omega(\bolda)} \\
        \vdots \\
        c_{d1} \boldt^{\boldp_0^{(1)}} s^{\omega(\boldp_0^{(1)})} + \cdots + 
        c_{dn} \boldt^{\boldp_0^{(n)}} s^{\omega(\boldp_0^{(n)})} -  
        c_{d0} s^{\omega(\boldzero)} 
        = \sum_{\bolda \in S} c_{d,\bolda} \boldt^{\bolda} s^{\omega(\bolda)} 
    \end{matrix}
    \right.
\end{equation}
which is constructed by multiplying each term in \eqref{equ:affine-cut-again}
by a rational power of the new variable $s$ whose exponent is determined by the
lifting function $\omega : S \to \mathbb{Q}$.
Clearly, $H(\boldt,1) = \boldzero$ is exactly the system \eqref{equ:affine-cut-again}
which we aim to solve (inside $(\Cstar)^d$).
As $s$ varies, however, $H$ represents a continuous deformation of  
the system \eqref{equ:affine-cut-again}, or a \emph{homotopy}.
The central idea behind the \emph{homotopy continuation method} 
for solving systems of equations is the deformation of a system into 
a ``starting system'' which one can solve easily.
Then numerical continuation methods are employed to trace the movement of 
the solutions of the starting system under the deformation toward the solutions
of the original system which one aims to solve.

The key here is to find an appropriate starting system that can be easily solved.
As is, $H(\boldt,0)$ cannot be used as the starting system since at $s=0$,
the system is either identically zero or undefined.
Therefore certain transformation is necessary to produce a meaningful and solvable
starting system.
Such transformations are given by the regular simplicial subdivision
discussed in \S \ref{sec:degree}.

Still let $\calD$ be a regular simplicial subdivision obtained by the algorithm
presented in \S \ref{sec:algo}.
Recall that each cell in $\calD$ is a projection of a cell of the form
$\{ \hat{\bolda}_0, \ldots, \hat{\bolda}_d \}$ such that 
$\conv \{ \hat{\bolda}_0, \ldots, \hat{\bolda}_d \}$ is a lower $d$-face
of $\conv \hatS$ that is characterized by \eqref{equ:lower-d-face}. 
That is, there exists a (unique) vector of the form 
$\hat{\boldalpha} = (\alpha_1,\ldots,\alpha_d,1)$ such that
\begin{equation} \label{equ:lower-d-face-again}
    \begin{aligned}
        \langle \hat{\bolda}_0, \hat{\boldalpha} \rangle &=
        \langle \hat{\bolda}_j, \hat{\boldalpha} \rangle & \text{ for } & j = 1, \dots, d \\
        \langle \hat{\bolda}_0, \hat{\boldalpha} \rangle &<
        \langle \hat{\bolda}  , \hat{\boldalpha} \rangle & \text{ for } & \text{all } \bolda \in S \;
    \end{aligned} .
\end{equation}
Using $\hat{\boldalpha} = (\alpha_1,\ldots,\alpha_d,1)$, 
we shall consider the change of variables 
\begin{equation} \label{equ:change-var}
    \boldt =
    \begin{cases}
        t_1 &= y_1 s^{\alpha_1} \\
            &\vdots \\
        t_d &= y_d s^{\alpha_d}
    \end{cases}
\end{equation}
with which $H$ becomes
\begin{equation*}
    H (\boldt,s) = H(y_1 s^{\alpha_1},\ldots, y_d s^{\alpha_d}, s) = 
    \left\{
    \begin{matrix}
        \sum_{\bolda \in S} c_{1,\bolda} \boldy^{\bolda} 
        s^{\langle \bolda, \boldalpha \rangle + \omega(\bolda)} =
        \sum_{\bolda \in S} c_{1,\bolda} \boldy^{\bolda} 
        s^{ \langle \hat{\bolda}, \hat{\boldalpha} \rangle } \\
        \vdots \\
        \sum_{\bolda \in S} c_{d,\bolda} \boldy^{\bolda} 
        s^{\langle \bolda, \boldalpha \rangle + \omega(\bolda)} =
        \sum_{\bolda \in S} c_{d,\bolda} \boldy^{\bolda} 
        s^{ \langle \hat{\bolda}, \hat{\boldalpha} \rangle } \\
    \end{matrix}
    \right.
\end{equation*}
Let $\beta = \langle \hat{\bolda}_0, \hat{\boldalpha} \rangle$ and define
a new homotopy
\begin{equation*}
    H^{\boldalpha,\beta} (\boldy,s) = 
    s^{-\beta} H(y_1 s^{\alpha_1},\ldots, y_d s^{\alpha_d}, s) =
    \left\{
    \begin{matrix}
        s^{-\beta} 
        \sum_{\bolda \in S} c_{1,\bolda} \boldy^{\bolda} 
        s^{ \langle \hat{\bolda}, \hat{\boldalpha} \rangle } \\
        \vdots \\
        s^{-\beta} 
        \sum_{\bolda \in S} c_{d,\bolda} \boldy^{\bolda} 
        s^{ \langle \hat{\bolda}, \hat{\boldalpha} \rangle } \\
    \end{matrix}
    \right.
\end{equation*}
Note that the new homotopy still has the necessary property that
$H^{\boldalpha,\beta} (\boldy,1) = \boldzero$ is identical to the system 
\eqref{equ:affine-cut-again} which we aim to solve.

One important observation here is that, by \eqref{equ:lower-d-face-again}, 
there are precisely $d+1$ terms in each component of $H^{\boldalpha,\beta}(\boldy,s)$ 
having no power of $s$ (the terms corresponding to $\bolda_0,\ldots,\bolda_d$), 
and all other terms have positive powers of $s$.
Consequently, at $s=0$, terms with positive powers of $s$ vanish, leaving only
\begin{equation} \label{equ:affine-cut-starting}
    \begin{cases}
        c_{1,\bolda_0} \boldy^{\bolda_0} +
        c_{1,\bolda_1} \boldy^{\bolda_1} + \cdots + 
        c_{1,\bolda_d} \boldy^{\bolda_d} & =  0 \\
        c_{2,\bolda_0} \boldy^{\bolda_0} +
        c_{2,\bolda_1} \boldy^{\bolda_1} + \cdots + 
        c_{2,\bolda_d} \boldy^{\bolda_d} & =  0 \\
        & \mathrel{\makebox[\widthof{=}]{\vdots}} \\
        c_{d,\bolda_0} \boldy^{\bolda_0} +
        c_{d,\bolda_1} \boldy^{\bolda_1} + \cdots + 
        c_{d,\bolda_d} \boldy^{\bolda_d} & =  0
    \end{cases}
\end{equation}
To simplify the notation, let
\begin{align*}
    C &=
    \begin{pmatrix}
    	c_{1,\bolda_0} & \cdots & c_{1,\bolda_d} \\
    	\vdots & \ddots & \vdots       \\
    	c_{d,\bolda_0} & \cdots & c_{d,\bolda_d}
    \end{pmatrix}
    &
    \Gamma &=
    \begin{pmatrix}
        \bolda_0 & \cdots & \bolda_d
    \end{pmatrix}
\end{align*}
then the above equation can be written as 
\begin{equation} \label{equ:affine-cut-starting-matrix}
    C \cdot (\boldy^\Gamma)^\top = \boldzero .
\end{equation}
For generic choices of the coefficients,
there exists a nonsingular matrix $G \in M_{d \times d}(\C)$ such that
\begin{equation}
    G C = 
    \begin{pmatrix}
    	c^*_{11} &          &        &          & c^*_{12} \\
    	         & c^*_{21} &        &          & c^*_{22} \\
    	         &          & \ddots &          & \vdots   \\
    	         &          &        & c^*_{d1} & c^*_{d2}
    \end{pmatrix} .
\end{equation}
for some $c^*_{ij} \in \Cstar$.
Then without altering its solution set, \eqref{equ:affine-cut-starting-matrix}
can be transformed into the equivalent system
\begin{equation} \label{equ:affine-cut-eliminated}
    G C (\boldt^\Gamma)^\top =
    \left\{
    \begin{matrix}
        c^*_{11} \boldy^{\bolda_0} & & & + & 
        c^*_{12} \boldy^{\bolda_d} & = & 0 \\
        & c^*_{21} \boldy^{\bolda_1} & & + &       
        c^*_{22} \boldy^{\bolda_d} & = & 0 \\
        & & & \mathrel{\makebox[\widthof{=}]{\vdots}} & \mathrel{\makebox[\widthof{=}]{\vdots}} & \mathrel{\makebox[\widthof{=}]{\vdots}} \\
        & & c^*_{d1} \boldy^{\bolda_{d-1}} & + &
        c^*_{d2} \boldy^{\bolda_d} & = & 0
    \end{matrix}
    \right.
\end{equation}
which is also a Laurent binomial system.
Therefore, the algorithm outlined can then be used to solve this system.
The solutions are precisely the solutions of the starting system \eqref{equ:affine-cut-starting}
for the homotopy $H^{\boldalpha,\beta}$.
Then numerical continuation techniques can be applied to trace the solutions
toward $s=1$ producing solutions to the target system \eqref{equ:affine-cut-again},
which will be points in the witness set of the component $V$ of $\V(\boldx^A-b)$.

Recall that the construction of the homotopy $H^{\boldalpha,\beta}$ depends on a cell
in the regular simplicial subdivision $\calD$ of $\conv S$.
It is typical for $\calD$ to contain more than one cell.
In this case, each cell induce a different homotopy of the form of $H^{\boldalpha,\beta}$.
The above construction is a special case of the \emph{polyhedral homotopy}
\cite{huber_polyhedral_1995}, and its theory guarantees that as one go through all the
cells in $\calD$, the resulting homotopies of the form $H^{\boldalpha,\beta}$
will find all the points in the witness set of $V$.


\section{Master space of $\mathcal{N}=1$ gauge theories} \label{sec:mspace}
In this section, we consider a system arising from theoretical physics, in particular, string theory.
A central area of current research in string theory is the study of the vacuum moduli space,
which, roughly speaking, is the space of continuous solutions (or the affine algebraic variety)
of a multivariate nonlinear function, called the superpotential of the theory under consideration. 
Here, the vacuum moduli spaces are spaces of special holonomy such as
Calabi-Yau or $G_2$ manifolds. 
Different positive-dimensional components of the vacuum moduli space 
correspond to different particle branches, such as mesonic, baryonic, etc. 
Symbolic algebraic geometry methods have been used to study the 
complicated structures of the vacuum moduli spaces of various 
string theory models \cite{Gray:2009fy,Gray:2008zs,Gray:2006gn}.
However, the methods are known to run out of the steam for even moderate sized systems 
due to the algorithmic complexity issues. 
Recently, numerical algebraic geometry methods have been introduced
to string theory research and have solved bigger systems 
\cite{Greene:2013ida,Hauenstein:2012xs,He:2013yk,MartinezPedrera:2012rs,Mehta:2009,Mehta:2011xs,Mehta:2011wj,Mehta:2012wk}. 

In this article, we consider special types of models coming from string theory 
in which the systems to be solved are binomial systems, and in which the vacuum 
moduli spaces are composed of unions of positive-dimensional components. 
We take a model which is actively investigated by string theorists 
because its vacuum moduli space is a combination of mesonic and baryonic branches 
\cite{Forcella:2008eh,Forcella:2008bb}.
Such moduli spaces are called \emph{master spaces}. 

In particular, we consider the superpotential for $\mathcal{N}=1$ gauge theories 
for a $D3$-brane on the Abelian orbifold $\mathbb{C}^3/\mathbb{Z}_{m} \times \mathbb{Z}_{k}$. 
The superpotential for this theory, for fixed $m,k \in \N$, is given by
\begin{equation} \label{equ:superpotential}
    W_{m,k} = \sum_{i =0}^{m-1} \sum_{j=0}^{k-1} 
    x_{i,j} y_{i+1,j} z_{i+1,j+1} - y_{i,j} x_{i,j+1} z_{i+1,j+1}
\end{equation}
where the periodic boundary conditions are imposed, e.g., 
$x_{i,m} = x_{i,0}$ for any $i$ and $x_{k,j} = x_{0,j}$ for any $j$. 
This is a polynomial in $3mk$ variables: $x_{i,j}, y_{i,j}, z_{i,j}$ 
for the combinations of $i \in \Z_k$ and $j \in \Z_m$.
For example, when $m=k=2$, the superpotential is
\begin{align*}
    W_{2,2} &= x_{0,0} y_{1,0} z_{1,1} - y_{0,0} x_{0,1} z_{1,1} 
             + x_{0,1} y_{1,1} z_{1,0} - y_{0,1} x_{0,0} z_{1,0} \\
            &+ x_{1,0} y_{0,0} z_{0,1} - y_{1,0} x_{1,1} z_{0,1}
             + x_{1,1} y_{0,1} z_{0,0} - y_{1,1} x_{1,0} z_{0,0},
\end{align*}
a polynomial in 12 variables 
$x_{0,0}$, $x_{0,1}$, $x_{1,0}$, $x_{1,1}$,
$y_{0,0}$, $y_{0,1}$, $y_{1,0}$, $y_{1,1}$,
$z_{0,0}$, $z_{0,1}$, $z_{1,0}$, $z_{1,1}$.
We are interested in finding the critical points of $W_{m,k}$,
that is, points at which all the partial derivatives of the superpotential $W_{m,k}$,
with respect to variables $x_{i,j},y_{i,j},z_{i,j}$, are zero.
These points are precisely the solutions to the system of polynomial equation
\begin{equation} \label{equ:mspace}
    \frac{\partial W_{m,k}}{\partial x_{i,j}} = \frac{\partial W_{m,k}}{\partial y_{i,j}} = \frac{\partial W_{m,k}}{\partial z_{i,j}} = 0
\end{equation}
in the variables $x_{i,j},y_{i,j},z_{i,j}$.

Notice that in $W_{m,k}$, each variable appears in exactly two distinct terms.
Consequently, the partial derivative of $W_{m,k}$ with respect to each variable consists of exactly two terms, 
hence it forms a binomial polynomial. For instance, 
\[
    \frac{\partial W_{2,2}}{\partial x_{0,0}} =  y_{1,0} z_{1,1} - y_{0,1} z_{1,0}, \quad
    \frac{\partial W_{2,2}}{\partial x_{0,1}} = -y_{0,0} z_{1,1} + y_{1,1} z_{1,0}. \quad
\]
Therefore \eqref{equ:mspace} is indeed a binomial system
which shall simply be denoted by $\nabla W_{m,k}$.
We are interested in computing the dimension and degree of components
of the $\Cstar$-solution set $\V(\nabla W_{m,k})$ of this system.

The dimension and the degree of the top dimensional components of this system
was first computed in \cite{Forcella:2008bb} for up to $m=3$ and $k=5$ using the 
Gr\"obner basis method. 
Later on, in \cite{Mehta:2012wk}, the dimensions and degrees of all the
components for up to $m=3=k$ were carried out using numerical algebraic geometry methods. 
The parallel GPU-based implementation of the binomial solver we have proposed can compute, 
very quickly, the dimension and the global parametrization of the $\Cstar$-solution set 
for higher values of $m$ and $k$. 
Table \ref{tab:mspace-dim} and \ref{tab:mspace-dim-long} show the dimension of 
$\V(\nabla W_{m,k}) \subset (\C^\ast)^{3mk}$ for a range of values for $m$ and $k$.
More importantly, our implementation shows impressive efficiency in computing the degree
of $\V(\nabla W_{m,k})$ for larger values of $m$ and $k$, including a component of degree
as high as $50467100$, for $m=4$ and $k=5$.
Table \ref{tab:mspace-deg} shows the degree of $\V(\nabla W_{m,k})$ for a range of
$m$ and $k$ values which is a significant expansion of the existing results
presented in \cite{Forcella:2008bb,Mehta:2012wk} and a substantial improvement over 
the existing algorithm outlined in \cite{chen_solutions_2014}.

\begin{table}
    \centering
    \begin{tabular}{|c||r|r|r|r|r|r|r|r|}
    	\hline
    	$m/k$ &   1 &  2 &  3 &  4 &  5 &  6 &  7 &  8 \\ \hline \hline
    	  1   & N/A &  4 &  5 &  6 &  7 &  8 &  9 & 10 \\ \hline
    	  2   &   4 &  6 &  8 & 10 & 12 & 14 & 16 & 18 \\ \hline
    	  3   &   5 &  8 & 11 & 14 & 17 & 20 & 23 & 26 \\ \hline
    	  4   &   6 & 10 & 14 & 18 & 22 & 26 & 30 & 34 \\ \hline
    	  5   &   7 & 12 & 17 & 22 & 27 & 32 & 37 & 42 \\ \hline
    	  6   &   8 & 14 & 20 & 26 & 32 & 38 & 44 & 50 \\ \hline
    	  7   &   9 & 16 & 23 & 30 & 37 & 44 & 51 & 58 \\ \hline
    	  8   &  10 & 18 & 26 & 34 & 42 & 50 & 58 & 66 \\ \hline
    \end{tabular}
    \vskip 1ex
    \caption{
        The dimension of $\V(\nabla W_{m,k})$ for a range of values for $m$ and $k$.
    }
    \label{tab:mspace-dim}
\end{table}

\begin{table}
    \centering
    \scriptsize
    \begin{tabular}{|c|r|r|r|r|r|r|r|r|r|r|r|r|r|r|r|r|}
    	\hline
    	$m=k$ &  9 &  10 &  11 &  12 &  13 &  14 &  15 &  20 &  25 &  30 &   35 &   40 \\ \hline\hline
    	Dim.  & 83 & 102 & 123 & 146 & 171 & 198 & 227 & 402 & 627 & 902 & 1227 & 1602 \\ \hline
    \end{tabular}
    \vskip 1ex
    \caption{
        The dimension of $\V(\nabla W_{m,k})$ for a range of larger values for $m=k$.
    }
    \label{tab:mspace-dim-long}
\end{table}
    
\begin{table}
    \centering
    \tiny
    \renewcommand{\tabcolsep}{5pt}
    \begin{tabular}{|c||c|c|c|c|c|c|c|c|}
    	\hline
    	$m/k$ &     1     &      2      &        3        &       4        &       5        &       6        &       7        &        8        \\ \hline\hline
    	  1   &           &  \known 2   &    \known 4     &    \known 8    &   \known 16    &   \known 32    &       64       &       128       \\ \hline
    	  2   & \known 2  &  \known 14  &    \known 92    &   \known 584   &  \known 3632   &     22304      &     135872     &     823424      \\ \hline
    	  3   & \known 4  &  \known 92  &   \known 1620   &  \known 26762  & \known 437038  &    7029180     &   111135118*   & $\ge 100100328$ \\ \hline
    	  4   & \known 8  & \known 584  &  \known 26762   &    1169876     &    50467100    & $\ge 11907022$ & $\ge 37567994$ &  \\ \hline
    	  5   & \known 16 & \known 3632 &  \known 437038  &    50467100    & $\ge 99710106$ & $\ge 62944504$ &                &  \\ \hline
    	  6   & \known 32 &    22304    &     7029180     & $\ge 11907022$ & $\ge 62944504$ &                &                &  \\ \hline
    	  7   &    64     &   135872    &   111135118*    & $\ge 37567994$ &                &                &                &  \\ \hline
    	  8   &    128    &   823424    & $\ge 100100328$ &                &                &                &                &  \\ \hline
    \end{tabular}
    \vskip 2ex
    \caption{
        The degree of the $\Cstar$-solution set defined by \eqref{equ:superpotential} 
        for a range of $m$ and $k$ values.
        This table lists only the results that can be computed within \emph{1 hour}
        on a \techy{NVidia GTX 780} graphics card with the double-precision version 
        of the GPU-based parallel algorithm presented in this work.
        Shaded entries correspond to and agree with the results already presented in 
        \cite{Forcella:2008bb,Mehta:2012wk}.
        Entries marked by * are results that cannot be computed with any CPU-based
        program within a reasonable amount of time (2 days for multi-core systems
        and 7 days for clusters).
        Entries marked by $\ge$ are lower bounds of the degrees computed by counting
        the total number of cells in the simplicial subdivision of the polytope
        associated with \eqref{equ:superpotential}.
    }
    \label{tab:mspace-deg}
\end{table}

Table \ref{tab:mspace-gpu-cpu} shows the speedup ratio achieved by the 
GPU-based algorithm, presented in \S \ref{sec:degree}, 
over its closest CPU-based implementation \techy{MixedVol-2.0} \cite{lee_mixed_2011}
which is widely regarded as one of fastest serial software program for computing 
``mixed volume'' (See Remark \ref{rmk:history} and Appendix \ref{sec:kushnirenko}
for its connection with degree computation considered in this article).
Remarkably, with sufficient GPU threads nearly 30 fold speedup ratio
has been achieved by the double-precision version of the algorithm.
When the single-precision version is used, even higher speedup ratio can be achieved,
Unfortunately, it appears that single-precision is, in general, 
not reliable in handling very large problems due to its insufficient precision.

\begin{table}
    \centering
    \begin{tabular}{|r|r|r|}
    	\hline
    	GPU threads & DP Speedup ratio & SP Speedup ratio \\ \hline\hline
    	         64 &             0.00 &             0.00 \\ \hline
    	        128 &             0.00 &             0.00 \\ \hline
    	        256 &             0.00 &             0.00 \\ \hline
    	        512 &             0.91 &             0.73 \\ \hline
    	       1024 &             0.98 &             4.14 \\ \hline
    	       2048 &             1.15 &             6.66 \\ \hline
    	       4096 &             2.20 &            10.99 \\ \hline
    	       8192 &             4.01 &            18.71 \\ \hline
    	   $2^{14}$ &             7.99 &            35.00 \\ \hline
    	   $2^{15}$ &            15.00 &            40.10 \\ \hline
    	   $2^{16}$ &            16.33 &            45.33 \\ \hline
    	   $2^{17}$ &            29.47 &            44.99 \\ \hline
    	   $2^{18}$ &            28.33 &            41.06 \\ \hline
    \end{tabular}
    \vskip 2ex
    \caption{
        Speedup ratios achieved by the GPU-based double-precision (DP) and single-precision (SP)
        algorithm respectively on \techy{NVidia GTX 780} when compared to \techy{MixedVol-2.0}.
        ``$0.00$'' represents speedup ratios too small to be measured reliably.
        The number of threads are chosen to be multiples of 32 which is the
        ``warp size'' (smallest group of threads in \techy{CUDA} framework).
    }
    \label{tab:mspace-gpu-cpu}
\end{table}

More important to note is the great potential of the GPU-based algorithm
when multiple GPU devices are used.
\autoref{tab:gpu-multi} shows the speedup ratio achieved by multiple GPU devices 
when compared to a single GPU, a single CPU, and a small cluster of 100 nodes.
With three GPU devices, over 60 fold speedup over the single-threaded CPU-based
algorithm (\techy{MixedVol-2.0}) has been achieved.
The most surprising result is the comparison between the GPU-based algorithm, 
developed in this article, running on three GPU devices and a similar CPU-based algorithm
running on a small cluster.
\techy{MixedVol-3} is a parallel version of \techy{MixedVol-2.0} \cite{lee_mixed_2011}
and, now, a part of a larger software program \techy{Hom4PS-3} \cite{chen_hom4ps-3:_2014}.
With three $\techy{NVidia GTX 780}$ our GPU-based algorithm computes the degree for
$\V(\nabla W_{4,5})$ faster than \techy{MixedVol-3} on a small cluster totaling 100 
\techy{Intel Xeon 2.4Ghz} processor cores.

In computing the degrees of $\V(\nabla W_{m,k})$ for certain larger $m$ and $k$,
while the GPU based algorithm was able to obtain the simplicial subdivisions
of the polytopes associated with $\V(\nabla W_{m,k})$, the volumes of certain 
cells, which are given as a matrix determinant  \eqref{equ:simplex-vol}, 
could not be computed with sufficient accuracy to ensure the exactness of
the answer.
However, since the volume of each cell is at least one,
the total number of cells is therefore a lower bound of the degree which equals
the total volume of all the cells.
Although these lower bounds are likely to be much smaller than the actual degrees,
given the sheer size of these systems, these partial results still merit further 
investigations and improvements on the approach presented here.
The lower bounds are therefore also included in Table \ref{tab:mspace-deg}
(entries marked with ``$\ge$'').

\begin{table}
    \centering
    \begin{tabular}{|r|r|r|r|}
    	\hline
    	              N.o. devices &     1 &     2 &     3 \\ \hline
    	Speedup over single device & 100\% & 188\% & 213\% \\ \hline
    	     Max. speedup over CPU & 28.33 & 54.12 & 61.00 \\ \hline
    	     Max. speedup over a cluster of 100 nodes & 0.48 & 0.91 & 1.04 \\ \hline
    \end{tabular}
    \vskip 1ex
    \caption{
        Speedup ratios achieved by using multiple identical \techy{NVidia GTX 780} devices
        with the single device performance (using the same algorithm) as a reference.
        }
    \label{tab:gpu-multi}
\end{table}

While the rigorous analysis and physical interpretation of the data presented here
are outside the scope of this article, the rich set of data shown in
Tables \ref{tab:mspace-dim}, \ref{tab:mspace-dim-long}, and \ref{tab:mspace-deg}
appear to show some general pattern.
To motivate further research in this important problem, we summarize these patterns
in the form of a conjecture:

\begin{conjecture}
    In general, for $m,k \in \Z^+$ with $m \ne k$, the solution set
    $\V(\nabla W_{m,k})$ consists of a single component of dimension
    \[
        \dim \V(\nabla W_{m,k}) = mk+2.
    \]
    Furthermore, for $m=1$ and $m=2$, the degree of the solution set is given by
    \begin{align*}
        \deg \V(\nabla W_{1,k}) &= 2 \cdot \deg \V(\nabla W_{1,k-1}) = 2^{k-1} \\
        \deg \V(\nabla W_{2,k}) &= 6 \cdot \deg \V(\nabla W_{2,k-1}) + 2^{2k-3} = 2 \cdot 6^{k-1} + \sum_{j=0}^{k-2} 2^{2(k-j)-3} \cdot 6^j 
    \end{align*}
\end{conjecture}

\section{Conclusion}

This this article, we proposed a parallel algorithm for computing the degree
of components of a $\Cstar$-solution set defined by a binomial system
that is specifically designed for GPU devices.
Numerical experiments with the \techy{CUDA} based implementation shows 
remarkable performance and scalability when applied to the binomial systems 
of the master space of $\mathcal{N}=1$ gauge theories.

\section*{Acknoweldgement}

DM was supported by a DARPA Young Faculty Award 
and an Australian Research Council DECRA fellowship. 
TC was supported in part by NSF under Grant DMS 11-15587.
TC and DM would like to thank Daniel Brake, Yang-Hui He and Thomas Kahle for their feedback on this paper.
TC would also like to thank Dirk Colbry for the helpful discussions
and the Institute for Cyber-Enabled Research at Michigan State University
for providing the necessary hardware and computational infrastructure.
\FloatBarrier

\begin{appendices}

\section{Kushnirenko's theorem} \label{sec:kushnirenko}

\begin{theorem}[\name{Kushnirenko} \cite{kushnirenko_newton_1975}]
    \label{thm:kushnirenko} 
    Consider the system of $k$ Laurent polynomial equations
    \[
        \begin{cases}
            c_{1,1} \boldx^{\bolda^{(1)}} + c_{1,2} \boldx^{\bolda^{(2)}} + \dots + c_{1,k} \boldx^{\bolda^{(\ell)}} &= 0 \\
            c_{2,1} \boldx^{\bolda^{(1)}} + c_{2,2} \boldx^{\bolda^{(2)}} + \dots + c_{2,k} \boldx^{\bolda^{(\ell)}} &= 0 \\
            &\vdots \\
            c_{k,1} \boldx^{\bolda^{(1)}} + c_{k,2} \boldx^{\bolda^{(2)}} + \dots + c_{k,k} \boldx^{\bolda^{(\ell)}} &= 0
        \end{cases}
    \]
    in $k$ variables $\boldx=(x_1,\dots,x_k)$ in which every equation has the same set 
    of monomials determined by exponent vectors $\bolda^{(1)},\dots,\bolda^{(\ell\,)} \in \Z^k$.
    With ``generic'' coefficients $c_{i,j} \in \C^\ast$,
    the solutions of this system in $(\C^\ast)^k$ are all isolated and nonsingular.
    The total number of these isolated solutions is
    \[
        k! \cdot \op{Vol}_k (\op{conv} \{ \bolda^{(1)}, \dots, \bolda^{(\ell\,)} \} ).
    \]
\end{theorem}

This important result was later generalized significantly in \cite{bernshtein_number_1975}
where the number of isolated nonzero $\Cstar$-solutions of a system of Laurent polynomial 
equation is shown to be equal to the \emph{mixed volume} of the \emph{Newton polytopes}
of the system.
This number is now commonly known as the \emph{BKK bound} of a Laurent polynomial system.
Therefore, the degree computation discussed above can be considered as a special case
of the BKK bound i.e. mixed volume computation.

\end{appendices}

\bibliography{bibliography_NPHC_NAG,MasterSpace}
\bibliographystyle{abbrv}

\end{document}